\pdfoutput=1
\RequirePackage{ifpdf}
\ifpdf 
\documentclass[pdftex]{sigma}
\else
\documentclass{sigma}
\fi

\usepackage{young}
\usepackage[enableskew,vcentermath]{youngtab}

\usepackage{mathrsfs}
\usepackage{amsxtra}
\usepackage{amscd}

\numberwithin{equation}{section}

\newtheorem{Theorem}{Theorem}[section]
\newtheorem{cor}[Theorem]{Corollary}
\newtheorem{lem}[Theorem]{Lemma}
\newtheorem{Claim}[Theorem]{Claim}
\newtheorem{prop}[Theorem]{Proposition}
\newtheorem{thm}[Theorem]{Theorem}

 { \theoremstyle{definition}
\newtheorem{defi}[Theorem]{Definition}
\newtheorem{Example}[Theorem]{Example}
\newtheorem{rem}[Theorem]{Remark} }

\newcommand{\Glie}{\mathfrak{g}} 
\newcommand{\Gaff}{\widehat{\mathfrak{g}}} 
\newcommand{\Gafft}{\widehat{\mathfrak{g}'}} 
\newcommand{\Uc}{U} 
\newcommand{\ev}{\mathrm{ev}} 
\newcommand{\CF}{\mathcal{F}} 
\newcommand{\CR}{\mathbf{R}} 
\newcommand{\BGG}{\mathcal{O}} 
\newcommand{\wCP}{\mathfrak{P}} 
\newcommand{\lCP}{\widehat{\mathfrak{P}}} 
\newcommand{\lCQ}{\widehat{\mathcal{Q}}}
\newcommand{\lwt}{\mathrm{wt}_{\ell}} 
\newcommand{\wt}{\mathrm{wt}} 
\newcommand{\Bw}{\mathbf{w}}
\newcommand{\Bn}{\mathbf{n}}
\newcommand{\Bm}{\mathbf{m}}
\newcommand{\Bf}{\mathbf{f}}
\newcommand{\CX}{\mathcal{X}}
\newcommand{\CY}{\mathcal{Y}}
\newcommand{\CEl}{\mathcal{E}_{\ell}} 
\newcommand{\BC}{\mathbb{C}} 
\newcommand{\BZ}{\mathbb{Z}} 
\newcommand{\BV}{\mathbf{V}} 
\newcommand{\BW}{\mathbf{W}}
\newcommand{\BQ}{\mathbf{Q}} 
\newcommand{\BP}{\mathbf{P}} 
\newcommand{\End}{\operatorname{End}}
\newcommand{\Id}{\operatorname{Id}} 
\newcommand{\Sm}{\mathbb{S}} 
\newcommand{\SF}{\mathscr{F}} 
\newcommand{\CG}{\mathscr{G}} 
\newcommand{\super}{\mathbb{Z}_2} 
\newcommand{\even}{\overline{0}} 
\newcommand{\odd}{\overline{1}} 
\newcommand{\CB}{\mathcal{B}} 
\newcommand{\hd}{\mathrm{hd}} 

\begin{document}

\allowdisplaybreaks

\newcommand{\arXivNumber}{1410.0837}

\renewcommand{\PaperNumber}{066}

\FirstPageHeading

\ShortArticleName{Asymptotic Representations of Quantum Af\/f\/ine Superalgebras}

\ArticleName{Asymptotic Representations\\ of Quantum Af\/f\/ine Superalgebras}

\Author{Huafeng ZHANG}

\AuthorNameForHeading{H.~Zhang}

\Address{Departement Mathematik and Institut f\"ur Theoretische Physik, ETH Z\"{u}rich, Switzerland}
\Email{\href{mailto:huafeng.zhang@math.ethz.ch}{huafeng.zhang@math.ethz.ch}}

\ArticleDates{Received April 21, 2017, in f\/inal form August 17, 2017; Published online August 22, 2017}

\vspace{-1mm}

\Abstract{We study representations of the quantum af\/f\/ine superalgebra associated with a~general linear Lie superalgebra. In the spirit of Hernandez--Jimbo, we construct inductive systems of Kirillov--Reshetikhin modules based on a cyclicity result that we established previously on tensor products of these modules, and realize their inductive limits as modules over its Borel subalgebra, the so-called $q$-Yangian. A new generic asymptotic limit of the same inductive systems is proposed, resulting in modules over the full quantum af\/f\/ine superalgebra. We derive generalized Baxter's relations in the sense of Frenkel--Hernandez for representations of the full quantum group.}


\Classification{17B37; 17B10; 81R50}

\vspace{-2.5mm}

\section{Introduction}
Let $q$ be a non-zero complex number which is not a root of unity. Let $\Glie := \mathfrak{gl}(M,N)$ be the {\it general linear Lie superalgebra}. Let $U_q(\Gaff)$ be the associated quantum af\/f\/ine superalgebra (Def\/inition~\ref{def: quantum affine superalgebra}). This is a Hopf superalgebra neither commutative nor co-commutative, and it can be seen as a deformation of the universal enveloping algebra of the af\/f\/ine Lie superalgebra $\Gaff := \Glie \otimes \BC[t,t^{-1}]$ of central charge zero without derivation.
In this paper we study a certain one-parameter family of inf\/inite-dimensional representations of $U_q(\Gaff)$, which arise from suitable limits of f\/inite-dimensional irreducible modules, the so-called {\it Kirillov--Reshetikhin modules}.

{\bf 1. Background.} Our study of quantum af\/f\/ine superalgebras is inspired, on the one hand, from the integrability structure of supersymmetric models (AdS/CFT, Hubbard model, etc.; see for example \cite{Beisert3}), and on the other hand, from the representation theoretical interpretation of transfer matrices for quantum integrable systems.

In the early seventies, towards the study of the transfer matrix $\mathbf{T}$ of the eight-vertex model, Baxter \cite{Baxter} introduced $\BQ$-operators and $\textbf{T-Q}$ functional relations to solve the spectra of $\mathbf{T}$. Within the framework of quantum inverse scattering method, the representation meaning of $\mathbf{Q}$-operators and $\textbf{T-Q}$ relations was clarif\/ied in a series of papers~\cite{BHK, BLZ1,BLZ2,BLZ} for lattice models and quantum f\/ield theories whose symmetry algebras are the quantum af\/f\/ine algebras attached to~$\mathfrak{sl}_2$,~$\mathfrak{sl}_3$. The idea goes roughly as follows. The af\/f\/ine quantum group $U_q(\widehat{\mathfrak{sl}_2})$ admits a~universal $R$-matrix $\mathcal{R} \in \mathfrak{B}_+ \otimes \mathfrak{B}_-$ with $\mathfrak{B}_{\pm}$ Borel subalgebras. One f\/ixes the so-called {\it quantum space}, a~representation $W$ of $\mathfrak{B}_-$ provided by the integrable models. This def\/ines an $\textbf{L}$-operator, an element of a completed tensor product $\mathfrak{B}_+ \otimes \End W$. The $\textbf{T},\BQ$-operators, as elements of $\End W$, are twisted traces of $\textbf{L}$ over various representations of $\mathfrak{B}_+$: f\/inite-dimensional evaluation representations over $U_q(\widehat{\mathfrak{sl}_2})$ for $\textbf{T}$ and {\it oscillator representations} for~$\BQ$. Baxter's $\textbf{T-Q}$ relations are then deduced from tensor product decompositions of representations of~$\mathfrak{B}_+$.

The oscillator representations of Borel subalgebras were subsequently extended to Lie superalgebras: $\mathfrak{sl}(2|1)$ in \cite{BT}; $\mathfrak{gl}(M,N)$ in \cite{Tsuboi}; twisted case of $\mathfrak{osp}(1|2)$ in \cite{IZ, KZ}. In these works, one of the main ingredients is the oscillator realization of Borel subalgebras.

The recent work \cite{BL} on hypergeometric equations asks for similar $\textbf{T-Q}$ formulations for the exceptional Lie superalgebra $D(2,1;\alpha)$. The associated quantum af\/f\/ine superalgebra is still less understood in aspects of structure/representation theory, beyond the Drinfeld loop realizations~\cite{HSTY}. Its universal $R$-matrix remains unknown.

In \cite{HJ}, Hernandez--Jimbo proposed the oscillator representations of Borel subalgebras for an arbitrary non-twisted quantum af\/f\/ine algebra $U_q(\widehat{\mathfrak{a}})$. Their main idea is to take a suitable inductive limit of a distinguished family of f\/inite-dimensional $U_q(\widehat{\mathfrak{a}})$-modules, the Kirillov--Reshetikhin modules. The limit construction enabled Frenkel--Hernandez~\cite{FH} to derive generalized $\textbf{T-Q}$ relations in terms of representations and to solve a conjecture of Frenkel--Reshetikhin on the spectra of quantum integrable systems~\cite{FR}.

In this paper, we extend Hernandez--Jimbo's limit construction to the quantum af\/f\/ine superalgebra $U_q(\Gaff)$. The Borel subalgebra in our situation will be the $q$-Yangian $Y_q(\Glie)$. Furthermore we perform a new limit, {\it generic asymptotic limit}, to the inductive system of Kirillov--Reshetikhin modules, resulting in modules over $U_q(\Gaff)$ itself.

{\bf 2. A toy example.} The generic asymptotic limit in the present paper is best viewed in the case of the f\/inite-type Drinfeld--Jimbo quantum group $U_q(\mathfrak{sl}_2)$. This is an associative algebra generated by three elements $e^{+}$, $e^-$, $K$ subject to relations
\begin{gather*} K\ \mathrm{invertible},\qquad K e^{\pm} = q^{\pm 2} e^{\pm}K,\qquad e^+e^- - e^-e^+ = \frac{K-K^{-1}}{q-q^{-1}}. \end{gather*}
To each positive integer $k$ is attached a representation of $U_q(\mathfrak{sl}_2)$ on $V_k := \oplus_{i=0}^k \BC v_i$ with
\begin{gather*} K v_i = q^{k-2i} v_i,\qquad e^+ v_i = \frac{q^i - q^{-i}}{q-q^{-1}} v_{i-1},\qquad e^- v_i = \frac{q^{k-i}-q^{-k+i}}{q-q^{-1}} v_{i+1}. \end{gather*}
The matrix entries of $K,e^{\pm} \in \operatorname{End}(V_k)$ at $(v_i,v_j)$ are Laurent polynomials in $q^k$ for $k > i+j+1$. Specializing such polynomials to a f\/ixed non-zero complex number $c \in \BC^{\times}$, we obtain a~representation of $U_q(\mathfrak{sl}_2)$ on $V_{\infty} := \oplus_{i=0}^{\infty} \BC v_i$ with\footnote{Let $V$ be a vector space with basis $B$ and let $b,b' \in B$. The matrix entry of a~linear endomorphism $f \in \operatorname{End}(V)$ at $(b',b)$ is by def\/inition the coef\/f\/icient of~$b'$ in~$f(b)$.}
\begin{gather*} K v_i = c q^{-2i} v_i,\qquad e^+ v_i = \frac{q^i - q^{-i}}{q-q^{-1}} v_{i-1},\qquad e^- v_i = \frac{cq^{-i}-c^{-1}q^{i}}{q-q^{-1}} v_{i+1}. \end{gather*}
Roughly speaking $V_{\infty}$ is an analytic continuation of the $(V_k)_{k\in \BZ_{>0}}$ with respect to the discrete parameter~$q^k$~-- one replaces~$q^k$ everywhere by~$c$.

{\bf 3. Main results.} Let $I_0 := \{1,2,\dots,M+N-1\}$ be the set of Dynkin vertices of $\Glie$. There are $U_q(\Gaff)$-valued power series $\phi_i^{\pm}(z)$ in $z^{\pm 1}$ for $i \in I_0$ whose coef\/f\/icients mutually commute; they can be viewed as $q$-analogs of $A \otimes t^{\pm n} \in \Gaff$
with $A$ being a diagonal matrix in $\Glie$ and $n$ a positive integer. The $q$-Yangian $Y_q(\Glie)$ contains the coef\/f\/icients of the $\phi_i^+(z)$, but not $\phi_i^-(z)$. There is a highest weight representation theory adapted to a triangular decomposition of $U_q(\Gaff)$ whose Cartan part is generated by the $\phi_i^{\pm}(z)$.

Fix $r$ a Dynkin vertex and $a \in \BC^{\times}$ a spectral parameter. To a positive integer $k$ is attached the Kirillov--Reshetikhin (KR) module, the unique f\/inite-dimensional irreducible $U_q(\Gaff)$-module~$W_k$ which is generated by a highest weight vector $\omega_k$ such that\footnote{In~\cite{HJ} actually $ \frac{q_r^k-zaq_r^{-k}}{1-za}$ was used in the inductive system of KR modules. This is due to opposite triangular properties of the coproduct of the $\phi_i^{\pm}(z)$ between \cite[Section~7]{Damiani} and Proposition~\ref{prop: JM elements}.\label{foot 1}}{\samepage
 \begin{gather*}\phi_r^{\pm}(z) \omega_k = q_r^k \frac{1-za}{1-zaq_r^{2k}} \omega_k,\qquad \phi_i^{\pm}(z) \omega_k = \omega_k \qquad \mathrm{for}\ \ i \neq r \end{gather*}
 as power series in $z^{\pm 1}$. Here $q_r = q$ for $r \leq M$ and $q_r = q^{-1}$ for $r > M$.}

The $W_k$ are af\/f\/ine analogs of the $U_q(\mathfrak{sl}_2)$-modules $V_k$ above. The uniform choice of bases for them relies on an inductive system of vector superspaces (not of $U_q(\Gaff)$-modules)
\begin{gather*}W_1 \subseteq W_2 \subseteq W_3 \subseteq \cdots,\qquad \omega_1 = \omega_2 = \omega_3 = \cdots. \end{gather*}
We follow the idea of Hernandez--Jimbo~\cite{HJ}, which is a fusion procedure relating dif\/ferent tensor products of KR modules. In the super case we need the main results in our previous paper~\cite{Z2} to validate the fusion (Lemma~\ref{lem: KR fund modules}).

Let $W_{\infty} = \cup_{k>0} W_k$ be the inductive limit. Choose a basis $B_k$ of $W_k$ inductively on $k > 0$ so that $\omega_1 \in B_1 \subseteq B_{2} \subseteq B_3 \subseteq \cdots$. Then $B_{\infty} := \cup_{k>0} B_k$ forms a basis of $W_{\infty}$.

{\bf 3.1. Asymptotic property.} To each triple $(x, b, b') \in U_q(\Gaff) \times B_{\infty} \times B_{\infty}$ is associated a~unique Laurent polynomial $P_{b'b}^x(u) \in \BC[u,u^{-1}]$ satisfying (Lemma~\ref{lem: asym property KR}):
\begin{itemize}\itemsep=0pt
\item the matrix entry of $x \in \operatorname{End}(W_k)$ at $(b',b)$ is $P_{b'b}^x\big(q_r^k\big)$ for $k$ large enough;\footnote{The $U_q(\Gaff)$-module structure on~$W_k$ def\/ines $x \in \operatorname{End}(W_k)$. When $k$ is large enough $b,b' \in B_k$ are basis vectors of $W_k$, so the matrix entry makes sense.}
\item given $b$ and $x$, we have $P_{b'b}^x(u) = 0$ for all but f\/initely many $b'$.
\end{itemize}
Our proof of the asymptotic property is more constructive than \cite{HJ}. It is based on representation theory of f\/inite-type quantum groups of $\mathfrak{sl}_2$ and $\mathfrak{sl}(1,1)$.

{\bf 3.2. Generic asymptotic limit.} Let $c$ be a non-zero complex number, called {\it spin para\-me\-ter}. The representation $(\rho_c,W_{\infty})$ of $U_q(\Gaff)$ is def\/ined by requiring (Corollary \ref{cor: generic asym repre}):
\begin{itemize}\itemsep=0pt
\item the matrix entry of $x \in \operatorname{End}(W_{\infty})$ at $(b',b)$ is the evaluation $P_{b'b}^x(c)$.
\end{itemize}
For example, let $b = \omega_1$ and $x = \phi_r^{\pm}(z)$. We have $P_{b'b}^x\big(q_r^k\big) = q_r^k \frac{1-za}{1-zaq_r^{2k}}\delta_{b'b}$ by def\/inition of KR modules. As a power series in $z^{\pm 1}$ its coef\/f\/icients are Laurent polynomials in $q_r^k$. Therefore $\phi_r^{\pm}(z) \omega_1 = c \frac{1-za}{1-zac^2} \omega_1$ in $(\rho_c, W_{\infty})$. Similarly $\phi_i^{\pm}(z) \omega_1 = \omega_1$ if $i \neq r$.

The generic asymptotic limit applied to non-twisted quantum af\/f\/ine algebras in the Appendix, we obtain $U_q(\widehat{\mathfrak{a}})$-modules $(\rho_c)_{c\in \BC^{\times}}$. They belong to the category $\BGG$ introduced by Hernandez~\cite{He} and further studied by Mukhin--Young~\cite{Mukhin}. If $c \notin \pm q^{\BZ}$, then $\rho_c$ is irreducible and is a~minimal af\/f\/inization of a parabolic Verma module~\cite{Mukhin}. Analytic continuation was used to prove genericity properties of minimal af\/f\/inizations~\cite{Mukhin}.

{\bf 3.3. Hernandez--Jimbo's limit.} The arguments of \cite{HJ} can be adapted to our situation to get a $Y_q(\Glie)$-module structure on $W_{\infty}$, an oscillator module in \cite{BT,Tsuboi}. One modif\/ies the $Y_q(\Glie)$-module structure on $W_k$, by tensoring with a~one-dimensional $Y_q(\Glie)$-module, so that \smash{$\phi_i^+(0) \omega_k = \omega_k$} for all $i \in I_0$. In this way one loses the $U_q(\Gaff)$-module structure.

 The advantage is that $P_{b'b}^x(u) \in \BC[u]$ for each triple $(x,b,b') \in Y_q(\Glie) \times B_{\infty} \times B_{\infty}$. We def\/ine the representation $(\rho_+, W_{\infty})$ of $Y_q(\Glie)$ by requiring (Corollary~\ref{cor: generic asym repre}):
 \begin{itemize}\itemsep=0pt
\item the matrix entry of $x \in \operatorname{End}(W_{\infty})$ at $(b',b)$ is the evaluation $P_{b'b}^x(0)$.
\end{itemize}
 Again take $x = \phi_r^+(z)$ and $b = \omega_1$. Then $P_{b'b}^x\big(q_r^k\big) = \frac{1-za}{1-zaq_r^{2k}} \delta_{b'b}$, which as a power series in $z$ has as coef\/f\/icients polynomials in $q_r^k$. It follows that $\phi_r^+(z) \omega_1 = (1-za) \omega_1$ in $(\rho_+, W_{\infty})$. Similarly $\phi_i^{\pm}(z) \omega_1 = \omega_1$ if $i \neq r$.

Informally, one can think of $(\rho_+, W_{\infty})$ as $(\rho_0,W_{\infty})$. Contrary to the non-graded case, $W_{\infty}$ is f\/inite-dimensional $\big(2^{MN}\big)$ for the odd Dynkin vertex $r = M$.

{\bf 3.4. Generalized Baxter's relations.} As in \cite{HJ} we introduce a monoidal category $\BGG$ of representations of $Y_q(\Glie)$ including all the f\/inite-dimensional $U_q(\Gaff)$-modules, the $\rho_+$ and $\rho_c$. In a~fraction ring of the Grothendieck ring of $\BGG$, the isomorphism class of a f\/inite-dimensional $U_q(\Gaff)$-module $[V]$ is a polynomial in the ratios $\frac{[(\rho_b,W_{\infty})]}{[(\rho_c,W_{\infty})]} $ whose coef\/f\/icients are isomorphism classes of one-dimensional $U_q(\Gaff)$-modules; see Theorem~\ref{thm: Baxter}. This form generalized Baxter's relations \`a~la Frenkel--Hernandez~\cite{FH}.

The key point in the proof is that normalized $q$-characters of the $(\rho_c, W_{\infty})$ are identical, which is a consequence of the generic asymptotic limit; see Lemma~\ref{lem: character asymptotic modules}. Note that the statement of Theorem~\ref{thm: Baxter} only involves $U_q(\Gaff)$-modules.

{\sloppy
{\bf 4.~Perspectives.} The generic asymptotic limit works for Felder's elliptic quantum groups~\cite{Felder} of which Borel subalgebras are still unknown; see~\cite{FZ,Z5}. It should eventually be done for other quantum af\/f\/ine superalgebras like $U_q\big(\widehat{D}(2,1;\alpha)\big)$.

}

In the sequels~\cite{FZ,Z4,Z5} we def\/ine Baxter $\BQ$-operators from transfer matrices of the $\rho_c$, and interpret Theorem~\ref{thm: Baxter} as generalized $\textbf{T-Q}$ relations of transfer matrices. Surprisingly the spin parameter $c$ becomes the spectral parameter of $\BQ$-operators.

This paper is organized as follows. Section~\ref{sec: idea} explains the idea of generic asymptotic construction. In Section~\ref{sec: 1}, we recall basic properties of the quantum af\/f\/ine superalgebra. Section~\ref{sec: KR} constructs inductive systems of Kirillov--Reshetikhin modules. In Section~\ref{sec: generic asymptotic representations} we carry out in detail the two limit constructions. In Section~\ref{sec: O and q} we introduce category $\BGG$ and state generalized Baxter's relations, whose proof is completed in Section~\ref{sec: GT bases} together with examples. In the appendix we apply our generic asymptotic construction to non-graded quantum af\/f\/ine algebras.

\section{Idea of asymptotic construction} \label{sec: idea}
Throughout this paper, all the vector spaces and algebras are def\/ined over the base f\/ield $\BC$. For two vector spaces $V$, $W$, let $\operatorname{Hom}(V,W)$ denote the set of all linear maps $V \longrightarrow W$.

Fix $A$ to be a unital associative algebra. Let $S$ be a system of algebraic generators, so that~$A$ is the quotient of the free associative algebra $\BC\langle S\rangle$ by the def\/ining ideal $\mathcal{I}_S$.

For all positive integer $k \in \BZ_{>0}$ let a representation $\rho^k\colon A \longrightarrow \operatorname{End}(V_k)$ of $A$ be given. Let $(F_{k,l}\colon V_l \longrightarrow V_k)_{l<k}$ be an inductive system of vector spaces; namely $F_{k,l} F_{l,m} = F_{k,m}$ for $m<l<k$ as linear maps $V_m \longrightarrow V_k$. Assume that the $F_{k,l}$ are {\it injective}.

Fix $L,K \in \BZ_{>0}$. We assume the {\it asymptotic property}: for all $l \in \BZ_{>0}$ and $s \in S$, there exists a $\operatorname{Hom}(V_l,V_{l+L})$-valued Laurent polynomial $P_{s;l}(u) = \sum\limits_{i=-K}^{K} P_{s;l}^{[i]} u^i$ in $u$ such that
\begin{gather*}\rho^k(s) F_{k,l} = F_{k,l+L} P_{s;l}(u)|_{u = q^k} \in \operatorname{Hom}(V_l,V_k) \qquad \mathrm{for}\quad k > l+L. \end{gather*}
 Let us prove that the Laurent polynomial $P_{s;l}(u)$ is unique. Indeed, let $Q(u)$ be another such Laurent polynomial and let $D(u) = P_{s;l}(u) - Q(u)$. Then from the injectivity of $F_{k,l+L}$ we have $D(u)|_{u=q^k} = 0$ for all $k > l+L$. Since $q$ is not a root of unity, a~Vandermonde matrix argument shows that $D(u) = 0$.

Similarly one shows that for f\/ixed $s \in S$ and $-K\leq i \leq K$, the linear maps $\big(P_{s;l}^{[i]}\big)_{l>0}$ form a morphism of inductive systems: $F_{l+L,m+L} P_{s;m}^{[i]} = P_{s;l}^{[i]} F_{l,m} \in \operatorname{Hom}(V_m, V_{l+L})$ for $m < l$. Let $P_s^{[i]} \in \End (V_{\infty})$ be its inductive limit with $V_{\infty}$ being the inductive limit of $(V_l, F_{k,l})$.

\begin{Claim} \label{Claim1}
Let $c \in \BC^{\times}$. Then $s \mapsto \sum\limits_{i=-K}^K P_s^{[i]} c^i$ defines a representation of $A$ on $V_{\infty}$.
\end{Claim}

The proof is again a Vandermonde matrix argument, and is omitted. As an example, suppose $s, t \in S$ and $s t \in \mathcal{I}_S$. Let us show
\begin{gather*}
\left(\sum_{i=-K}^K P_s^{[i]} c^i\right)\left( \sum_{i=-K}^K P_t^{[i]} c^i\right) = 0.
\end{gather*} We have
\begin{gather*}
0 = \rho^k(st) F_{k,l} = \rho^k(s) \underline{\rho^k(t) F_{k,l}} = \underline{\rho^k(s) F_{k,l+L}} P_{t;l}\big(q^k\big) = F_{k,l+2L} P_{s;l+L}\big(q^k\big)P_{t;l}\big(q^k\big)
\end{gather*}
for all $k > l+2L$. This forces $P_{s;l+L}(u)P_{t;l}(u)|_{u = q^k} = 0$ for all $k > l+2L$ and so $P_{s;l+L}(u) P_{t;l}(u)$ $= 0$. Taking inductive limit $l \rightarrow \infty$ and $u = c$ leads to the desired identity.

If furthermore the $P_{s;l}(u)$ are polynomials in $u$, then in Claim~\ref{Claim1} one can take $c = 0$.
\section{Backgrounds on quantum superalgebras} \label{sec: 1}
This section collects basic facts on the RTT realization of the quantum af\/f\/ine superal\-ge\-bra~$U_q(\Gaff)$, the $q$-Yangian $Y_q(\Glie)$ and the quantum superalgebra $\Uc_q(\Glie)$ following \cite{Z2}.

 Fix $M,N \in \BZ_{>0}$. Set $I := \{1,2,\dots,M+N\},\ \super := \BZ/2\BZ = \{\even,\odd\}$. For $i \in I$,
\begin{gather} \label{def: shift}
 |i| =: \begin{cases}
\even, & i\leq M, \\
\odd, & i > M,
\end{cases} \qquad d_i := \begin{cases}
\hphantom{-}1, & i \leq M,\\
-1, & i > M,
\end{cases} \qquad q_i := q^{d_i}.
\end{gather}
Def\/ine the weight lattice $\BP := \oplus_{i \in I} \BZ \epsilon_i$ with bilinear form $(\,,\,)\colon \BP \times \BP \longrightarrow \BZ, (\epsilon_i,\epsilon_j) = \delta_{ij}d_i$. Let $|\cdot|\colon \BP \longrightarrow \super$ be the morphism of abelian groups such that $|\epsilon_i| = |i|$. Set $I_0 := I \setminus \{M+N\}$. For $i \in I_0$, let $\alpha_i := \epsilon_i - \epsilon_{i+1}$. Def\/ine the root lattice $\BQ = \oplus_{i\in I_0} \BZ \alpha_i \subset \BP$, and root cones $\BQ_{+} := \oplus_{i \in I_0} \BZ_{\geq 0} \alpha_i$ and $\BQ_- := - \BQ_+$.

{\it Only} three cases of $|x| \in \super$ will be admitted: $x \in I$; $x \in \BP$; $x$ is a $\super$-homogeneous vector of a vector superspace $V$. Naturally $\operatorname{Hom}(V,V) =: \End(V)$ is a superalgebra.

Let $\BV := \bigoplus_{i \in I} \BC v_i$ be the vector superspace with parity $|v_i| = |i|$. The superalgebra $\End (\BV)$ has a basis formed of elementary matrices $E_{ij}\colon v_k \mapsto \delta_{jk} v_i$. Note that $|E_{ij}| = |i| + |j|$ and $E_{ij} E_{kl} = \delta_{jk} E_{il}$.
 Recall the {\it Perk--Schultz matrix} from \cite{Perk-Schultz}:
\begin{gather*}
R(z,w) := \sum\limits_{i\in I}\big(zq_i - wq_i^{-1}\big) E_{ii} \otimes E_{ii} + (z-w) \sum\limits_{i \neq j} E_{ii} \otimes E_{jj} \\
\hphantom{R(z,w) :=}{} + z \sum\limits_{i<j} \big(q_i-q_i^{-1}\big) E_{ji} \otimes E_{ij} + w \sum\limits_{i<j}\big(q_j-q_j^{-1}\big) E_{ij} \otimes E_{ji}.
\end{gather*}

\begin{defi}\cite{Z2} \label{def: quantum affine superalgebra} The quantum af\/f\/ine superalgebra $U_q(\Gaff)$ is def\/ined by
\begin{itemize}\itemsep=0pt
\item[(R1)] generators $s_{ij}^{(n)}$, $t_{ij}^{(n)}$ for $i,j \in I$ and $n \in \BZ_{\geq 0}$;
\item[(R2)] parity $\big|s_{ij}^{(n)}\big| = \big|t_{ij}^{(n)}\big| = |i| + |j|$;
\item[(R3)] RTT-relations \cite{FRT2, FRT} in $U_q(\Gaff) \otimes \End (\BV^{\otimes 2})[[z,z^{-1},w,w^{-1}]]$:
\begin{gather*}
R_{23}(z,w) T_{12}(z) T_{13}(w) = T_{13}(w) T_{12}(z) R_{23}(z,w), \\
R_{23}(z,w) S_{12}(z) S_{13}(w) = S_{13}(w) S_{12}(z) R_{23}(z,w), \\
R_{23}(z,w) T_{12}(z) S_{13}(w) = S_{13}(w) T_{12}(z) R_{23}(z,w), \\
t_{ij}^{(0)} = s_{ji}^{(0)} = 0 \qquad \mathrm{for}\quad 1 \leq i < j \leq M+N, \\
t_{ii}^{(0)} s_{ii}^{(0)} = 1 = s_{ii}^{(0)} t_{ii}^{(0)} \qquad \mathrm{for}\quad i \in I.
\end{gather*}
\end{itemize}
Here {\samepage
\begin{gather*}
T(z) = \sum_{i,j \in I} t_{ij}(z) \otimes E_{ij} \in (U_q(\Gaff) \otimes \End \BV)\big[\big[z^{-1}\big]\big], \\ t_{ij}(z) = \sum_{n \in \BZ_{\geq 0}} t_{ij}^{(n)} z^{-n} \in U_q(\Gaff)\big[\big[z^{-1}\big]\big]
\end{gather*} (similar convention for $S(z)$ with the $z^{-n}$ replaced by the $z^{n}$).}

The $q$-Yangian $Y_q(\Glie)$ is the subalgebra of $U_q(\Gaff)$ generated by the $s_{ij}^{(n)}$, $\big(s_{ii}^{(0)}\big)^{-1}$.

The quantum supergroup $U_q(\Glie)$ is the subalgebra of $U_q(\Gaff)$ generated by the $s_{ij}^{(0)}$ and $t_{ij}^{(0)}$. We write $s_{ij}^{(0)}$, $t_{ij}^{(0)}$ as $s_{ij}$, $t_{ij}$ when no confusion with the series $s_{ij}(z)$, $t_{ij}(z)$ arises.
\end{defi}

Usual convention: if $A$, $B$, $C$ are superalgebras and $T = \sum_i a_i \otimes b_i \in A \otimes B$, then we write $T_{12} := \sum_i a_i \otimes b_i \otimes 1 \in A \otimes B \otimes C$, $T_{13} := \sum_{i} a_i \otimes 1 \otimes b_i \in A \otimes C \otimes B$ and $T_{23} := \sum_i 1 \otimes a_i \otimes b_i \in C \otimes A \otimes B$.

$U_q(\Gaff)$ has a Hopf superalgebra structure with counit $\varepsilon\colon U_q(\Gaff) \longrightarrow \BC$ def\/ined by $\varepsilon\big(s_{ij}^{(n)}\big) = \varepsilon\big(t_{ij}^{(n)}\big) = \delta_{ij}\delta_{n0}$, and coproduct $\Delta\colon U_q(\Gaff) \longrightarrow U_q(\Gaff)^{\otimes 2}$:
\begin{gather*} 
\Delta \big(s_{ij}^{(n)}\big) = \sum_{m=0}^n \sum_{k \in I} \epsilon_{ijk} s_{ik}^{(m)} \otimes s_{kj}^{(n-m)}, \qquad \Delta \big(t_{ij}^{(n)}\big) = \sum_{m=0}^n \sum_{k \in I} \epsilon_{ijk} t_{ik}^{(m)} \otimes t_{kj}^{(n-m)}.
\end{gather*}
Here $\epsilon_{ijk} := (-1)^{|E_{ik}||E_{kj}|}$. The antipode $\Sm\colon U_q(\Gaff) \longrightarrow U_q(\Gaff)$ is determined by
\begin{eqnarray*}
&&(\Sm \otimes \operatorname{Id})(S(z)) = S(z)^{-1}, \qquad (\Sm \otimes \operatorname{Id})(T(z)) = T(z)^{-1}.
\end{eqnarray*}
 Notice that $Y_q(\Glie)$ and $U_q(\Glie)$ are sub-Hopf-superalgebras of $U_q(\Gaff)$.

 We need the Drinfeld--Jimbo generators $e_i^{\pm}, K_i \in U_q(\Glie)$ for $i \in I_0$ \cite[Proposition 3.3]{Z2}:
\begin{gather} \label{def: Drinfeld-Jimbo generator}
e_i^+ := \frac{s_{ii}^{-1}s_{i,i+1}}{1-q_i^{-2}},\qquad e_i^- := \frac{t_{i+1,i}t_{ii}^{-1}}{1-q_i^2},\qquad K_i := s_{ii}s_{i+1,i+1}^{-1}.
\end{gather}

Let us recall the relations of these generators from \cite{Z2}:\footnote{We use the super bracket $[x,y] := x y - (-1)^{|x||y|} y x$. Equation~\eqref{rel: relations between e_0 and f_i} is \cite[equation~(A.20)]{Z2} applied to by an evaluation map in Proposition~\ref{prop: evaluation morphism} below. $\Glie'$ is not to be confused with the derived algebra of~$\Glie$. }
\begin{gather}
 [e_i^+, e_j^- ] = \delta_{ij}\frac{K_i-K_i^{-1}}{q_i-q_i^{-1}} \qquad \mathrm{for}\quad i,j \in I_0, \label{rel: relations between e_i and f_j} \\
 [t_{ji},t_{kj}] = \big(q_j-q_j^{-1}\big)t_{jj}t_{ki},\qquad [s_{ij},s_{jk}] = \big(q_j-q_j^{-1}\big) s_{jj}s_{ik} \qquad \mathrm{if}\quad i < j < k, \label{rel: relations between e_0 and f_i} \\
 s_{ii}^{(0)} s_{jk}^{(n)} = q^{(\epsilon_i,\epsilon_j - \epsilon_k)} s_{jk}^{(n)}s_{ii}^{(0)},\qquad s_{ii}^{(0)}t_{jk}^{(n)} = q^{(\epsilon_i,\epsilon_j - \epsilon_k)} t_{jk}^{(n)}s_{ii}^{(0)}. \label{rel: Cartan relations for weight grading}
\end{gather}
Set $\Glie' := \mathfrak{gl}(N,M)$. Let us def\/ine the quantum af\/f\/ine superalgebra $U_q(\Gafft)$ in the same way as $U_q(\Gaff)$, except that $M$, $N$ are interchanged. Let $s_{ij}'^{(n)}$, $t_{ij}'^{(n)}$ for $i,j \in I$ and $n \in \BZ_{\geq 0}$ be the corresponding RTT generators of $U_q(\Gafft)$, so that their parities are $\big|s_{ij}'^{(n)}\big| = \big|t_{ij}'^{(n)}\big| = |i|' + |j|'$ where $|i|' = \even$ for $1 \leq i \leq N$ and $\odd$ otherwise. For $i,j \in I$, set
\begin{gather*}
\varepsilon_{ij} := (-1)^{|i|(|i|+|j|)},\qquad \varepsilon_{ij}' := (-1)^{|i|'(|i|'+|j|')}, \qquad \widehat{i} := M+N+1-i.
\end{gather*}
\begin{prop} \label{prop: evaluation morphism}
The following assignments define morphisms of superalgebras:
\begin{gather*}
\ev_a\colon \ U_q(\Gaff) \longrightarrow \Uc_q(\Glie),\qquad s_{ij}(z) \mapsto s_{ij} - za t_{ij},\qquad t_{ij}(z) \mapsto t_{ij} - z^{-1}a^{-1}s_{ij}; \\
\phi_{(f(z),g(z))}\colon \ U_q(\Gaff) \longrightarrow U_q(\Gaff),\qquad s_{ij}(z) \mapsto f(z) s_{ij}(z),\qquad t_{ij}(z) \mapsto g(z) t_{ij}(z); \\
\Phi_a\colon \ U_q(\Gaff) \longrightarrow U_q(\Gaff),\qquad s_{ij}^{(n)} \mapsto a^n s_{ij}^{(n)},\qquad t_{ij}^{(n)} \mapsto a^{-n} t_{ij}^{(n)}; \\
\Psi\colon \ U_q(\Gaff) \longrightarrow U_q(\Gaff)^{\mathrm{cop}},\qquad s_{ij}^{(n)} \mapsto \varepsilon_{ji}t_{ji}^{(n)},\qquad t_{ij}^{(n)} \mapsto \varepsilon_{ji}s_{ji}^{(n)}; \\
 \CF\colon \ U_q(\Gafft) \longrightarrow U_q(\Gaff)^{\mathrm{cop}},\qquad s_{ij}'^{(n)} \mapsto \varepsilon_{ji}' s_{\widehat{j}\widehat{i}}^{(n)},\quad t_{ij}'^{(n)} \mapsto \varepsilon_{ji}' t_{\widehat{j}\widehat{i}}^{(n)}.
\end{gather*}
Here $a \in \BC^{\times}$ and $f(z),g(z^{-1}) \in 1 + z\BC[[z]]$. The last three maps are Hopf superalgebra isomorphisms. For $(A,\Delta,\varepsilon)$ a Hopf superalgebra, $(A^{\mathrm{cop}},\Delta^{\mathrm{cop}},\varepsilon)$ denotes another Hopf superalgebra with the same underlying superalgebra $A$ but with twisted coproduct $\Delta^{\mathrm{cop}} := c_{A,A} \Delta$, where $c_{A,A}\colon A^{\otimes 2} \longrightarrow A^{\otimes 2}$ is the graded permutation $x \otimes y \mapsto (-1)^{|x||y|} y \otimes x$.
\end{prop}

 $\ev_a$ is called an {\it evaluation map} as $\ev \circ \iota = \Id_{\Uc_q(\Glie)}$. The maps $\phi_{(f(z),g(z))}$, $\Phi_a$, $\CF$ restrict to $q$-Yangians. Relation \eqref{rel: Cartan relations for weight grading} gives rise to the weight grading on $U_q(\Gaff)$: for $\alpha \in \BQ$,
\begin{gather*}
U_q(\Gaff)_{\alpha} = \big\{ x \in U_q(\Gaff) \, |\, s_{ii}^{(0)} x \big(s_{ii}^{(0)}\big)^{-1} = q^{(\epsilon_i,\alpha)} x \ \textrm{for}\ i \in I \big\}.
\end{gather*}
For example $s_{ij}^{(n)}$, $t_{ij}^{(n)}$ are of weight $ \epsilon_i - \epsilon_j$. This induces weight gradings on $Y_q(\Glie)$ and $U_q(\Glie)$.

We end this section with some facts on representations of $U_q(\Glie)$ following \cite{BKK}.

A $\Uc_q(\Glie)$-module $V$ admits a {\it weight grading} if it is a direct sum of weight spaces
\begin{gather} \label{def: weight space}
V_{\alpha} := \big\{ x \in V \, |\, s_{ii} x = q^{(\epsilon_i,\alpha)} x \ \textrm{for}\ i \in I \big\} \qquad \mathrm{with}\quad \alpha \in \BP.
\end{gather}
By equation~\eqref{rel: Cartan relations for weight grading}, $\Uc_q(\Glie)_{\alpha} V_{\beta} \subseteq V_{\alpha + \beta}$ for $\alpha,\beta \in \BP$. If furthermore all the weight spaces $V_{\alpha}$ are f\/inite-dimensional, then the {\it character} of $V$ can be def\/ined:
\begin{gather*}
\chi (V) := \sum_{\alpha \in \BP} \dim (V)_{\alpha} [\alpha] \in \BZ^{\BP}.
\end{gather*}
Here $\BZ^{\BP}$ is the abelian group of functions $\BP \longrightarrow \BZ$ and $[\alpha]\colon \beta \mapsto \delta_{\alpha,\beta}$.

 Let $\lambda \in \BP$. Up to isomorphism, there exists a unique irreducible $\Uc_q(\Glie)$-module, denoted by~$L(\lambda)$, which is generated by a vector $v_{\lambda}$ satisfying
\begin{gather*}
|v_{\lambda}| = |\lambda|,\!\!\qquad s_{kk} v_{\lambda} = q^{(\epsilon_k,\lambda)} v_{\lambda}, \!\!\qquad t_{kk} v_{\lambda} = q^{-(\epsilon_k,\lambda)} v_{\lambda},\!\!\qquad s_{ij} v_{\lambda} = 0, \!\!\qquad i,j,k \in I,\quad i < j.
\end{gather*}
$L(\lambda)$ is weight graded with f\/inite-dimensional weight spaces. $L(\lambda)_{\lambda} = \BC v_{\lambda}$ and $L(\lambda)_{\alpha} \neq 0$ only if $\lambda-\alpha \in \BQ_+$. (The proof, parallel to the non-graded case, is based on the triangular decomposition and PBW basis of $U_q(\Glie)$ in~\cite{Yam}.)

The f\/irst example is the {\it vector representation} $\pi$ of $\Uc_q(\Glie)$ on $\BV$ \cite[Example~1]{Z2}:
\begin{gather}
\pi(s_{ii}) = q_i E_{ii} + \sum_{j \neq i} E_{jj} =\pi \big(t_{ii}^{-1}\big) \qquad \textrm{for}\quad i \in I, \nonumber\\
\pi(s_{ij}) = \big(q_i - q_i^{-1}\big) E_{ij},\qquad \pi(t_{ji}) = \big(q_i^{-1}-q_i\big) E_{ji}\qquad \textrm{for}\quad 1 \leq i < j \leq M+N.\label{def: vector repr}
\end{gather}
We have $(\pi,\BV) \cong L(\epsilon_1)$ with $v_1 = v_{\epsilon_1}$ and $\chi (L(\epsilon_1)) = \sum\limits_{i \in I} [\epsilon_i]$.

\begin{defi} \cite{Z2} \label{def: KR}
 Kirillov--Reshetikhin module $W_{k,a}^{(r)}$ for $r \in I_0$, $k \in \BZ_{>0}$, $a \in \BC^{\times}$ is
\begin{itemize}\itemsep=0pt
\item[(1)] either $\ev_{aq^{2k}}^* L(k\varpi_r)$ with $\varpi_r = \sum\limits_{j=1}^r \epsilon_j$ and $r \leq M$,
\item[(2)] or $\ev_a^* (L(k \varpi_r) ) \otimes \BC_{|k \varpi_r|} $ with $\varpi_r = -\sum\limits_{j=r+1}^{M+N} \epsilon_j$ and $r > M$.
\end{itemize}
Here $\BC_{s}$ for $s \in \super$ is the one-dimensional $U_q(\Gaff)$-module $\varepsilon\colon U_q(\Gaff) \longrightarrow \BC$ of parity~$s$.
\end{defi}

Let $r,k \in \BZ_{>0}$ be such that $r \leq M$ or $k \leq N$. Consider the rectangle Young diag\-ram~$Y_{r,k}$ with~$r$ rows and~$k$ columns. We view $Y_{r,k}$ as a subset of $\BZ_{>0} \times \BZ_{>0}$ so that~$(i,j) \in Y_{r,k}$ corresponds to the box at $i$-th row and $j$-th column. An~{\it $(M,N)$-hook semi-standard tableau} of shape~$Y_{r,k}$ is a~function $T\colon Y_{r,k} \longrightarrow I = \{1 < 2 < \cdots < M+N\}$ such that:
\begin{itemize}\itemsep=0pt
\item[(i)] the entries in each row and column are weakly increasing;
\item[(ii)] the entries in $\{1,2,\dots,M\}$ are strictly increasing in each column;
\item[(iii)] the entries in $\{M+1,M+2,\dots, M+N\}$ are strictly increasing in each row;
\end{itemize}
Let $\CB_{r,k}$ be the set of all such functions.

\begin{thm}[\cite{BKK}] \label{thm: BKK Schur-Weyl duality}
For $1\leq r \leq M$ and $k \in \BZ_{>0}$
\begin{gather} \label{equ: character of simple modules}
\chi (L(k\varpi_r)) = \sum_{T \in \CB_{r,k}} \bigg[\sum_{(i,j) \in Y_{r,k}} \epsilon_{T(i,j)}\bigg].
\end{gather}
If $r = M$ and $k \geq N$, then $\dim L(k\varpi_M) = 2^{MN}$.
\end{thm}
As a consequence, $L(k\varpi_r)$ is f\/inite-dimensional for $1\leq r \leq M$. This is also true for $M < r \leq M+N$. Indeed, the pullback of the $U_q(\Glie)$-module $L(k\varpi_r)$ by $\CF\colon U_q(\Glie') \longrightarrow U_q(\Glie)$ in Proposition \ref{prop: evaluation morphism} is an irreducible module over $U_q(\Glie')$ of highest weight $k \varpi_{M+N-r}'$ (we add prime to distinguish $\Glie' = \mathfrak{gl}(N,M)$ with $\Glie$) so that its character can be computed by equation~\eqref{equ: character of simple modules} in terms of $(N,M)$-hook semi-standard tableaux of shape $Y_{M+N-r,k}$.

In \cite{BKK}, the tableaux correspond to Kashiwara's crystal basis of $L(\lambda)$. As an example: $\Glie = \mathfrak{gl}(2,2)$ and $\lambda = 2\epsilon_1+2\epsilon_2$, the tableaux are: $\young(11,22)$, $\young(11,23)$, $\young(11,24)$, $\young(11,34)$, $\young(12,23)$, $\young(12,24)$, $\young(12,34)$, $\young(13,23)$, $\young(13,24)$, $\young(13,34)$, $\young(14,24)$, $\young(14,34)$, $\young(22,34)$, $\young(23,34)$, $\young(24,34)$, $\young(34,34)$.

\section{Inductive system of Kirillov--Reshetikhin modules} \label{sec: KR}
We construct inductive system of KR modules, based on our previous result~\cite{Z2} and establish its asymptotic property as $U_q(\Glie)$-modules.

Let $V$ be a $U_q(\Gaff)$-module. A $\super$-homogeneous non-zero vector $v \in V$ is called a {\it highest $\ell$-weight vector} if it is a common eigenvector for the $s_{ii}^{(n)},t_{ii}^{(n)}$ and it is annihilated by the~$s_{ij}^{(n)}$,~$t_{ij}^{(n)}$ with $i < j$. $V$ is called a {\it highest $\ell$-weight module} if $V = U_q(\Gaff)v$ for some highest $\ell$-weight vector~$v$. In this case, $v$~is unique up to scalar multiple, and $V$ admits a unique irreducible quotient, called the {\it head} of $V$ and denoted by~$\hd(V)$.

For example, let $v_{\lambda} \in L(\lambda)$ be as in Section~\ref{sec: 1}, and let $a \in \BC^{\times}$. The evaluation module $\ev_a^*L(\lambda)$ contains a highest $\ell$-weight vector $w := \ev_a^*(v_{\lambda})$ with
\begin{gather*}
|w| = |\lambda|,\qquad s_{ii}(z) w = \big(q^{(\epsilon_i,\lambda)} - z a q^{-(\epsilon_i,\lambda)}\big) w,\qquad t_{ii}(z) w = \big(q^{-(\epsilon_i,\lambda)} - z^{-1}a^{-1} q^{(\epsilon_i,\lambda)}\big) w.
\end{gather*}

The tensor product of two highest $\ell$-weight vectors is also a highest $\ell$-weight vector. Let~$V$,~$V'$ be $U_q(\Gaff)$-modules. We write $V \simeq V'$ if there exists a one-dimensional $U_q(\Gaff)$-module $D$ such that $V \cong V' \otimes D$ as $U_q(\Gaff)$-modules. In this case, since $s_{ij}^{(n)} |_D = t_{ij}^{(n)} |_D = 0$ for $i \neq j$, we have that $V$ is of highest $\ell$-weight if and only if so is~$V'$.

\begin{thm}[{\cite[Theorem 5.2]{Z2}}] \label{thm: cyclicity of fund modules}
Let $r \in I_0$ and $a_1,a_2,\dots,a_k \in \BC^{\times}$. The $U_q(\Gaff)$-module $W_{1,a_1}^{(r)} \otimes W_{1,a_2}^{(r)} \otimes \cdots \otimes W_{1,a_k}^{(r)}$ is of highest $\ell$-weight if $\frac{a_i}{a_j} \notin q_r^{\BZ_{<0}}$ for $1\leq i < j \leq k$.\footnote{This result has been made stronger in the author's later works~\cite{Z3,Z4}.}
\end{thm}

If a sequence $(X_j)_{j\in \BZ}$ of $Y_q(\Glie)$-modules is given and $m, n \in \BZ$ with $m \leq n$, then we write
\begin{gather*} X_n \otimes X_{n-1} \otimes \cdots \otimes X_{m+1} \otimes X_m =: \bigotimes\limits^{\leftarrow}_{m\leq j \leq n} X_j. \end{gather*}
For example, $\bigotimes\limits^{\leftarrow}_{0\leq j \leq k-1} W_{1,aq_r^{2j}}^{(r)}$ is of highest $\ell$-weight for $k \in \BZ_{>0}$.
\begin{lem} \label{lem: KR fund modules}
Let $r \in I_0, a \in \BC^{\times}$ and $k \in \BZ_{>0}$. Then $\hd\Big(\bigotimes\limits^{\leftarrow}_{0\leq j \leq k-1} W_{1,aq_r^{2j}}^{(r)}\Big) \simeq W_{k,a}^{(r)}$.
\end{lem}
\begin{proof}
Let us compare the highest $\ell$-weight vectors $v,v'$ at the left-hand and right-hand sides respectively. Write $s_{ii}(z) v = f_i(z) v$ and $s_{ii}(z) v' = f_i'(z) v'$ for $i \in I$. By Def\/inition~\ref{def: KR} and Proposition~\ref{prop: evaluation morphism}, $f_i(z) = f_{i+1}(z)$ and $f_i'(z) = f_{i+1}'(z)$ for $i \in I_0 \setminus \{r\}$,
\begin{gather*} \frac{f_r(z)}{f_{r+1}(z)} = \prod_{j=1}^k \frac{1-zaq_r^{2k-2j}}{q_r^{-1}-zaq_r^{2k-2j+1}} = \frac{1-za}{q_r^{-k}-zaq_r^{k}} = \frac{f_r'(z)}{f_{r+1}'(z)}. \end{gather*}
Similar statement holds for eigenvalues of the $t_{ii}(z)$. Set $g^{\pm} = \prod\limits_{j=1}^{k-1}\big(1-z^{\pm 1} a^{\pm 1} q_r^{\pm 2j}\big)$. Then $\hd\Big(\bigotimes\limits^{\leftarrow}_{0\leq j \leq k-1} W_{1,aq_r^{2j}}^{(r)}\Big) \cong W_{k,a}^{(r)} \otimes \phi_{(g^+,g^-)}^*(\BC_{\even})$; see Def\/inition~\ref{def: KR}(2).
\end{proof}

From now on up to the end of Section~\ref{sec: generic asymptotic representations}, $r \in I_0$ and $a \in \BC^{\times}$ are f\/ixed.

For $k \in \BZ_{>0}$, let $\rho^k$ denote the representation of $U_q(\Gaff)$ on $W_{k,a}^{(r)}$ and let us f\/ix a highest $\ell$-weight vector $\omega_k$ in $W_{k,a}^{(r)}$. If $l,k \in \BZ_{>0}$ and $l < k$, then def\/ine the $U_q(\Gaff)$-module
\begin{gather*}
Z_{kl} := \phi^*_{((1-zaq_r^{2l})^{-1},(1-z^{-1}a^{-1}q_r^{-2l})^{-1})} W_{k-l,aq_r^{2l}}^{(r)},
\end{gather*}
and f\/ix a highest $\ell$-weight vector $\omega_{kl}$ in $Z_{kl}$. Let us show that $Z_{kl} \otimes W_{l,a}^{(r)}$ is of highest $\ell$-weight. By Lemma \ref{lem: KR fund modules}, one may replace $Z_{kl}$ and $W_{l,a}^{(r)}$ by the heads of the tensor products
$T_1 := \bigotimes\limits^{\leftarrow}_{l \leq j\leq k-1} W_{1,aq_r^{2j}}^{(r)}$ and $T_2 := \bigotimes\limits^{\leftarrow}_{0 \leq j\leq l-1} W_{1,aq_r^{2j}}^{(r)}$ respectively. By Theorem \ref{thm: cyclicity of fund modules}, $T_1 \otimes T_2$ is of highest $\ell$-weight, so is $\hd(T_1) \otimes \hd(T_2)$ as its quotient.

From the proof of Lemma \ref{lem: KR fund modules} follows $\hd\big(Z_{kl} \otimes W_{l,a}^{(r)}\big) \cong W_{k,a}^{(r)}$. As in \cite[Section~4]{HJ}, let $\SF_{k,l}\colon Z_{kl} \otimes W_{l,a}^{(r)} \longrightarrow W_{k,a}^{(r)}$ be the quotient map sending $\omega_{kl} \otimes \omega_l$ to $\omega_k$, and def\/ine the restriction map $F_{k,l}\colon W_{l,a}^{(r)} \longrightarrow W_{k,a}^{(r)}, x \mapsto \SF_{k,l}(\omega_{kl} \otimes x)$. In particular $F_{k,l}(\omega_l) = \omega_k$.

\begin{prop} \label{prop: first properties of inductive system}
 The maps $\big(F_{k,l}\colon W_{l,a}^{(r)} \longrightarrow W_{k,a}^{(r)}\big)_{l<k}$ verify the following properties.
\begin{itemize}\itemsep=0pt
\item[$(1)$] $F_{k,l}\big(W_{l,a}^{(r)}\big)_{l\varpi_r - \beta} \subseteq \big(W_{k,a}^{(r)}\big)_{k\varpi_r - \beta}$ for $\beta \in \BQ_+$, and $\rho^k(e_i^+) F_{k,l} = F_{k,l} \rho^l(e_i^+)$ for $i \in I_0$.
\item[$(2)$] $F_{k,l}\colon W_{l,a}^{(r)} \longrightarrow W_{k,a}^{(r)}$ is injective, and $F_{k,l} F_{l,m} = F_{k,m}$ for $m<l<k$.
\item[$(3)$] $\rho^k(e_i^-) F_{k,l} = F_{k,l} \rho^l(e_i^-)$ for $i \in I_0 \setminus \{r\}$. For $m > 0$, there exist linear maps $A_m, B_m\colon W_{m,a}^{(r)}$ $\longrightarrow W_{m+1,a}^{(r)}$ of parity $|\alpha_r|$ such that
\begin{gather*}
\rho^k(e_r^-) F_{k,m} = F_{k,m+1} \big(q_r^kA_m+ q_r^{-k} B_m\big) \colon \ W_{m,a}^{(r)} \longrightarrow W_{k,a}^{(r)} \qquad \mathrm{for}\quad k > m+1.
\end{gather*}
\end{itemize}
\end{prop}
\begin{proof}
(1) By def\/inition $\omega_{kl}$ is of weight $(k-l)\varpi_r$. Since the $U_q(\Gaff)$-linear map $\SF_{k,l}$ respects the weight gradings, $F_{k,l}$ changes the weights by $(k-l)\varpi_r$. By equation~\eqref{def: Drinfeld-Jimbo generator},
\begin{gather*}
\Delta (e_i^+) = 1 \otimes e_i^+ + e_i^+ \otimes K_i^{-1},\qquad \Delta(e_i^-) = K_i \otimes e_i^- + e_i^- \otimes 1.
\end{gather*}
Since $e_i^+ \omega_{kl} = 0$, we have: for $x \in W_{l,a}^{(r)}$,
\begin{gather*} e_i^+ F_{k,l}(x) = e_i^+ \SF_{k,l}(\omega_{kl} \otimes x) = \SF_{k,l}(\Delta(e_i^+)(\omega_{kl}\otimes x)) = \SF_{k,l}(\omega_{kl} \otimes e_i^+x) = F_{k,l}e_i^+(x). \end{gather*}

(2) Assume $\ker(F_{k,l}) \neq 0$. By (1) there exists $\mu \in \BP$ such that $\ker(F_{k,l})_{\mu + \alpha_i} = 0$ and $\ker(F_{k,l})_{\mu} \neq 0$ for $i \in I_0$. This implies $e_i^+ \ker(F_{k,l})_{\mu} = 0$ for all $i \in I_0$. By equation~\eqref{rel: relations between e_0 and f_i}, $s_{ij} \ker(F_{k,l})_{\mu} = 0$ for $1\leq i < j \leq M+N$. Since $W_{l,a}^{(r)} \cong L(l\varpi_r) \otimes \BC_{|l\varpi_r|}$ is an irreducible $U_q(\Glie)$-module, we have $\ker(F_{k,l})_{\mu} \subseteq (W_{l,a}^{(r)})_{l\varpi_r} = \BC \omega_l$, in contradiction with $F_{k,l} (\omega_l) = \omega_k$.

Consider the $U_q(\Gaff)$-module $S:= Z_{kl} \otimes Z_{lm} \otimes W_{m,a}^{(r)}$. It contains a highest $\ell$-weight vector $\omega := \omega_{kl} \otimes \omega_{lm} \otimes \omega_m$. By Lemma~\ref{lem: KR fund modules} we have as $U_q(\Gaff)$-modules
\begin{gather*} S \simeq \bigg(\bigotimes\limits^{\leftarrow}_{l\leq j\leq k-1} W_{1,aq_r^{2j}}^{(r)}\bigg) \otimes \bigg(\bigotimes\limits^{\leftarrow}_{m\leq j\leq l-1} W_{1,aq_r^{2j}}^{(r)}\bigg) \otimes \bigg(\bigotimes\limits^{\leftarrow}_{0\leq j\leq m-1}W_{1,aq_r^{2j}}^{(r)}\bigg). \end{gather*}
Theorem \ref{thm: cyclicity of fund modules} applied to the right-hand side, $S$ and $Z_{kl} \otimes Z_{lm}$ are of highest $\ell$-weight with heads~$W_{k,a}^{(r)}$ and~$Z_{km}$ respectively. Let $\CG_{k,m}^l\colon Z_{kl} \otimes Z_{lm} \longrightarrow Z_{km}$ be the quotient map sending~$\omega_{kl} \otimes \omega_{lm}$ to~$\omega_{km}$. We obtain $U_q(\Gaff)$-module morphisms from $S$ to $W_{k,a}^{(r)}$:
\begin{gather*} F := \SF_{k,l}(\operatorname{Id}_{Z_{kl}} \otimes \SF_{l,m}), \qquad G:= \SF_{k,m}\big(\CG_{k,m}^l \otimes \operatorname{Id}_{W_{l,a}^{(r)}}\big). \end{gather*}
Since $F(\omega) = G(\omega) = \omega_k$ and since $S$ is generated by $\omega$, we have $F = G$. Applying $F$, $G$ to $\omega_{kl} \otimes \omega_{lm} \otimes W_{m,a}^{(r)}$ gives the desired identity of (2).

(3) Assume $i \in I_0 \setminus \{r\}$. Since $e_i^+ \omega_l = 0$, $K_i \omega_l = \omega_l$, by~\eqref{rel: relations between e_i and f_j}, $e_i^+ e_i^- \omega_l = 0$. If $j \in I_0 \setminus \{i\}$, then $e_j^+ e_i^- \omega_l \in W_{l,a}^{(r)}$ is of weight $l\varpi_r - \alpha_i + \alpha_j \notin l \varpi_r + \BQ_-$ and is zero. So $e_i^- \omega_l$ is annihilated by all the $e_j^+$. As in~(2), $e_i^- \omega_l = 0$. Then $e_i^- F_{k,l} = F_{k,l} e_i^-$ as in~(1), using $\Delta(e_i^-)$.

For $i = r$, we adapt the proof of \cite[Lemma 7.6]{Z4}. We f\/ix $l = m+1$ in (2). Applying $F = G$ to $\omega_{kl} \otimes w' \otimes w \in S$ for $w' \in Z_{lm}$ and $w \in W_{m,a}^{(r)}$ and $k > l$ gives the identity:
\begin{gather} \label{(*)} F_{k,l}\SF_{l,m}(w' \otimes w) = \SF_{k,m}(\CG_{k,m}^l(\omega_{kl} \otimes w') \otimes w). \end{gather}
We compute $e_r^- \omega_{km}$ via the projection $\CG_{km}^l(\omega_{kl} \otimes \omega_{lm}) = \omega_{km}$ and $K_r \omega_{kl} = q_r^{k-l} \omega_{kl}$
\begin{gather*} e_r^- \omega_{km} = \CG_{k,m}^l\big(q_r^{k-l}\omega_{kl} \otimes e_r^- \omega_{lm} + e_r^- \omega_{kl} \otimes \omega_{lm}\big). \end{gather*}
Next consider the following vector in $Z_{kl} \otimes Z_{lm}$ of weight $(k-m)\varpi_r - \alpha_r$:
\begin{gather*}x := q_r^{-1}\frac{q_r^{k-l}-q_r^{l-k}}{q_r-q_r^{-1}} \omega_{kl} \otimes e_r^- \omega_{lm} - e_r^- \omega_{kl} \otimes \omega_{lm}. \end{gather*}
Based on $\Delta(e_r^+)$ and $l = m+1$ one checks that $e_r^+ x = 0$. If $j \in I_0 \setminus \{r\}$, then $e_j^+ x$ is of weight $(k-m)\varpi_r - \alpha_r + \alpha_j \notin (k-m)\varpi_r + \BQ_-$ and is zero. So $\CG_{k,m}^l(x) \in Z_{km}$ is annihilated by the $e_j^+$ for $j \in I_0$. As in (2), $\CG_{k,m}^l(x) = 0$. It follows that
\begin{gather*} e_r^- \omega_{km} = \CG_{k,m}^l\big(q_r^{k-l}\omega_{kl} \otimes e_r^- \omega_{lm} + e_r^- \omega_{kl} \otimes \omega_{lm} + x\big) = \frac{q_r^{k-m}-q_r^{m-k}}{q_r-q_r^{-1}} \CG_{k,m}^l(\omega_{kl} \otimes e_r^- \omega_{lm}). \end{gather*}
Taking $w' := e_r^- \omega_{lm} \in Z_{lm}$ in \eqref{(*)}, we compute
\begin{gather*}
e_r^- F_{k,m}(w)= e_r^- \SF_{k,m}(\omega_{km} \otimes w) = \SF_{k,m}(e_r^- \omega_{km} \otimes w) + q_r^{k-m} \SF_{k,m}(\omega_{km} \otimes e_r^- w) \\
\hphantom{e_r^- F_{k,m}(w)}{} = \frac{q_r^{k-m}-q_r^{m-k}}{q_r-q_r^{-1}} \SF_{k,m}\big(\CG_{k,m}^l(\omega_{kl} \otimes e_r^- \omega_{lm}) \otimes w\big) + q_r^{k-m}F_{k,m}(e_r^- w) \\
\hphantom{e_r^- F_{k,m}(w)}{}= \frac{q_r^{k-m}-q_r^{m-k}}{q_r-q_r^{-1}} F_{k,l} \SF_{l,m}( e_r^- \omega_{lm} \otimes w) + q_r^{k-m} F_{k,m}(e_r^- w).
\end{gather*}
This shows that $\rho^k(e_r^-) F_{k,m} = F_{k,m+1} \big(q_r^k A_m + q_r^{-k} B_m\big)$ where
\begin{gather*}
A_m(w) := q_r^{-m} F_{m+1,m}(e_r^-w) + \frac{q_r^{-m}}{q_r-q_r^{-1}} \SF_{m+1,m}(e_r^-\omega_{m+1,m}\otimes w), \\
B_m(w) := \frac{-q_r^{m}}{q_r-q_r^{-1}} \SF_{m+1,m}(e_r^-\omega_{m+1,m}\otimes w),
\end{gather*}
as linear maps $W_{m,a}^{(r)} \longrightarrow W_{m+1,a}^{(r)}$ are clearly independent of $k$ and of parity~$|\alpha_r|$.
\end{proof}

\begin{rem} \label{rem: asym finite} By equation~\eqref{rel: relations between e_0 and f_i}, for $1\leq i < j \leq M+N$, $s_{ii}^{-1}s_{ij}$ (resp. $t_{ji}t_{ii}^{-1}$) is a sum of monomials of the $e_{i}^+, e_{i+1}^+,\dots, e_{j-1}^+$ (resp.\ the $e_i^-, e_{i+1}^-,\dots, e_{j-1}^-$) where each $e_{h}^{\pm}$ for $i \leq h < j$ appears once. It follows that $\rho^k\big(s_{ii}^{-1}s_{ij}\big) F_{k,l} = F_{k,l} \rho^l\big(s_{ii}^{-1}s_{ij}\big)$ and:
\begin{itemize}\itemsep=0pt
\item[(i)] if $i > r$ or $j \leq r$, then $\rho^k\big(t_{ji}t_{ii}^{-1}\big) F_{k,l} = F_{k,l} \rho^l\big(t_{ji}t_{ii}^{-1}\big)$ for $k > l$;
\item[(ii)] if $i\leq r < j$, then $\rho^k\big(t_{ji}t_{ii}^{-1}\big) F_{k,l} = F_{k,l+1}\big(q_r^k A_{ji;l} + q_r^{-k} B_{ji;l}\big)$ for $k > l+1$, where~$A_{ji;l}$,~$B_{ji;l}$ are linear maps $W_{l,a}^{(r)} \longrightarrow W_{l+1,a}^{(r)}$ of parity $|\epsilon_j-\epsilon_i|$.
\end{itemize}
For example, $A_{r+1,r;l} = A_l$ and $B_{r+1,r;l} = B_l$. Furthermore, for $i \in I$ and $l < k$ we have
\begin{gather*} \rho^k(s_{ii}) F_{k,l} = q^{(\epsilon_i,(k-l)\varpi_r)} F_{k,l} \rho^l(s_{ii}),\qquad \rho^k(t_{ii}) F_{k,l} = q^{-(\epsilon_i,(k-l)\varpi_r)} F_{k,l} \rho^l(t_{ii}). \end{gather*}
\end{rem}

\section{Asymptotic representations} \label{sec: generic asymptotic representations}
We apply the asymptotic constructions of Section~\ref{sec: idea} ($c \in \BC^{\times}$ or $c = 0$) to the inductive system of KR modules in Proposition~\ref{prop: first properties of inductive system}, and obtain $U_q(\Gaff)$ or $Y_q(\Glie)$-modules. Set $T := \big\{s_{ij}^{(n)},t_{ij}^{(n)}\big\}$ (resp.~$S := \big\{s_{ij}^{(n)}\big\}$) to be the system of algebraic generators of~$U_q(\Gaff)$ (resp.~$Y_q(\Glie)$).

For $k > 0$, the map $S \ni s_{ij}^{(n)} \mapsto q^{(\epsilon_i,-k\varpi_r)} \rho^k\big(s_{ij}^{(n)}\big) \in \End\big(W_{k,a}^{(r)}\big)$ extends uniquely to a representation $\big(\tilde{\rho}^k, W_{k,a}^{(r)}\big)$ of~$Y_q(\Glie)$. Indeed, there is a representation $\big(\theta_k^{(r)},\BC\big)$ of $Y_q(\Glie)$ on the one-dimensional vector superspace of even parity given by $\theta_k^{(r)}\big(s_{ij}^{(n)}\big) = q^{(\epsilon_i,-k\varpi_r)} \delta_{ij} \delta_{n,0}$ and $(\tilde{\rho}^k, W_{k,a}^{(r)}) \cong (\theta_k^{(r)},\BC) \otimes \big(\rho^k,W_{k,a}^{(r)}\big)$ as $Y_q(\Glie)$-modules.\footnote{The $\tilde{\rho}^k$ are normalized Kirillov--Reshetikhin modules in~\cite{FH}.}
\begin{lem}\label{lem: asym property KR}
Let $t \in T,s \in S$ and $l > 0$. There exist $\operatorname{Hom}\big(W_{l,a}^{(r)},W_{l+1,a}^{(r)}\big)$-valued Laurent polynomials $P_{t;l}(u) = \sum\limits_{i=-2}^2 P_{t;l}^{[i]} u^i$ and $Q_{s;l}(u) = \sum\limits_{i=0}^1 Q_{s;l}^{[i]} u^i$ in $u$ of parity $|t|$ and $|s|$ respectively such that for all $k > l+1$ the following identities hold in $\operatorname{Hom}\big(W_{l,a}^{(r)}, W_{k,a}^{(r)}\big)$:
\begin{gather*}
 \rho^k(t) F_{k,l} = F_{k,l+1} P_{t;l}(u)|_{u = q_r^k},\qquad \tilde{\rho}^k(s) F_{k,l} = F_{k,l+1} Q_{s;l}(u)|_{u = q_r^{2k}}.
\end{gather*}
\end{lem}
\begin{proof}
If $r \leq M$, then for $i,j \in I$ we have as power series in $z$:
\begin{gather*}
 \sum_{n\geq 0} z^n \rho^k\big(s_{ij}^{(n)}\big) F_{k,l} = \rho^k(s_{ij}) F_{k,l} - zaq^{2k} \rho^k(t_{ij}) F_{k,l} \\
\hphantom{\sum_{n\geq 0} z^n \rho^k\big(s_{ij}^{(n)}\big) F_{k,l}}{} = q^{(k-l)(\epsilon_i,\varpi_r)} F_{k,l} \rho^l(s_{ij}) - z aq^{2k-(k-l)(\epsilon_j,\varpi_r)} \rho^k\big(t_{ij}t_{jj}^{-1}\big) F_{k,l} \rho^l(t_{jj}),
\end{gather*}
and $\rho^k\big(t_{ij}^{(n)}\big) = - a^{-1}q^{-2k} \rho^k\big(s_{ij}^{(1-n)}\big)$. If $r > M$, then $\rho^k\big(t_{ij}^{(n)}\big) = - a^{-1}\rho^k\big(s_{ij}^{(1-n)}\big)$ and
\begin{gather*}
\sum_{n\geq 0} z^n \rho^k\big(s_{ij}^{(n)}\big) F_{k,l} = q^{(k-l)(\epsilon_i,\varpi_r)}\rho^k(s_{ij}) F_{k,l} - za q^{-(k-l)(\epsilon_j,\varpi_r)} \rho^k\big(t_{ij}t_{jj}^{-1}\big) F_{k,l} \rho^l(t_{jj}).
\end{gather*}
We compute the right-hand sides in four cases based on Remark~\ref{rem: asym finite}:
\begin{gather} \label{tab: asym formula}
\begin{tabular}{|c|c|c|}
 \hline
 $r \leq M$ & $r > M$ & \tsep{1pt}\bsep{1pt}\\
 \hline
 $q^{k}\times F_{k,l} \rho^l\big(q^{-l}s_{ij}\big)$ & $F_{k,l} \rho^l(s_{ij})$ & $i \leq r < j$\tsep{1pt}\bsep{1pt} \\
 \hline
 $F_{k,l}\rho^l\big(s_{ij}-zaq^{2k}t_{ij}\big)$ & $ q^k \times F_{k,l} \rho^l\big(q^{-l}s_{ij} - z a q^{-2k+l}t_{ij}\big)$ & $r < i,j$ \tsep{1pt}\bsep{1pt}\\
 \hline
 $-za F_{k,l+1} \big(q^{2k+l}A_{ij;l} + q^{l} B_{ij;l}\big)$ & $ q^k \times (-za) F_{k,l+1} \big(q^{-2k}A_{ij;l} + B_{ij;l}\big)$ & $j \leq r < i$ \tsep{1pt}\bsep{1pt}\\
 \hline
 $q^k \times F_{k,l} \rho^l\big(q^{-l}s_{ij} - zaq^{l} t_{ij}\big)$ & $F_{k,l} \rho^l(s_{ij} - za t_{ij})$ & $i,j \leq r$ \tsep{1pt}\bsep{1pt}\\
 \hline
\end{tabular}
\end{gather}
Removing the factors ``$q^k \times$'' gives $\tilde{\rho}^k(s_{ij}(z)) F_{k,l}$, a polynomial in $q_r^{2k}$ of degree $\leq 1$.
\end{proof}

Let $\big(F_l\colon W_{l,a}^{(r)} \longrightarrow W_{\infty}\big)_{l>0}$ be the inductive limit of $\big(F_{k,l}\colon W_{l,a}^{(r)} \longrightarrow W_{k,a}^{(r)}\big)$.
\begin{cor} \label{cor: generic asym repre}Fix $c \in \BC^{\times}$. There exist representations $\rho_c$ of $U_q(\Gaff)$ and $\rho_+$ of $Y_q(\Glie)$ on $W_{\infty}$ defined by the formulas: for $t \in T, s \in S$,
\begin{gather*}\rho_c(t) = \lim\limits_{l\rightarrow\infty} P_{t;l}(u)|_{u = c},\qquad \rho_+(s) = \lim\limits_{l\rightarrow\infty} Q_{s;l}(u)|_{u = 0}. \end{gather*}
\end{cor}
 Let $W_{a;c}^{(r)}$ and $L_{r,a}^+$ denote the $U_q(\Gaff)$-module $(\rho_c,W_{\infty})$ and $Y_q(\Glie)$-module $(\rho_+,W_{\infty})$ respectively. For computational purpose the next observation is useful.\footnote{In~\cite{Z4}, Drinfeld second realization arising from a dif\/ferent Gauss decomposition from Section~\ref{sec: O and q} is used to resolve the issue in footnote~\ref{foot 1}. This results in dif\/ferent parameterizations of highest $\ell$-weights. The two asymptotic limits of KR modules are denoted by $\mathcal{W}_{c,a}^{(r)}, L_{r,a}^-$ therein and match with \cite{HJ}.}

\begin{rem} \label{rem: asym computation}
Let us be in the situation of Section~\ref{sec: idea}. For $x \in A$, one can f\/ind $K(x), L(x) \in \BZ_{>0}$ (depending on $x$ and $L, K$) such that: for all $l > 0$ there exists a $\operatorname{Hom}(V_l, V_{l+L(x)})$-valued Laurent polynomial $\sum\limits_{i=-K(x)}^{K(x)} P_{x;l}^{[i]} u^i = P_{x;l}(u)$ in $u$ with
\begin{gather*}\rho^k(x) F_{k,l} = F_{k,l+L(x)} P_{x;l}(u)|_{u = q^k} \in \operatorname{Hom}(V_l,V_k) \qquad \mathrm{for}\quad k > l+L(x). \end{gather*}
In the representation of $A$ on $V_{\infty}$, $x$ acts as $\lim\limits_{l \rightarrow \infty} P_{x;l}(c)$.
For example, let $x = s t$ with $s,t \in S$. Then $L(x) = 2L, K(x) = 2K$ and $P_{st;l}(u) = P_{s;l+L}(u) P_{t;l}(u)$.

To compute $x v$ with $v \in V_{\infty}$, one f\/inds $l>0$ such that $v = F_l(v')$ for some $v' \in V_l$; here $F_l\colon V_{l} \longrightarrow V_{\infty}$ is a structural map of the inductive limit. Then one writes $\rho^k(x) F_{k,l}(v') = F_{k,l+L(x)} \mathcal{P}_l(u)|_{u= q^k}$ for $k > l+L(x)$, where $\mathcal{P}_l(u)$ is a $V_{l+L(x)}$-valued Laurent polynomial in $u$ of degree bounded by $K(x)$. At last, $F_{l+L(x)}(\mathcal{P}_{l}(c))$ is exactly $x v \in V_{\infty}$.
\end{rem}

Consider the action of the $s_{ii}$ on $W_{a;c}^{(r)}$. Let $v = F_l(v')$ with $v' \in \big(W_{l,a}^{(r)}\big)_{l\varpi_r-\beta}$ and $\beta \in \BQ_+$. We have $F_{k,l}(v') \in \big(W_{k,a}^{(r)}\big)_{k\varpi_r-\beta}$ and so $\rho^k(s_{ii}) F_{k,l}(v') = q^{(\epsilon_i, k\varpi_r - \beta)} F_{k,l}(v')$. This gives
\begin{gather*}
\rho_c(s_{ii}) v = c^{d_r (\epsilon_i,\varpi_r)} q^{(\epsilon_i,-\beta)} v,\quad \rho_+(s_{ii}) v = q^{(\epsilon_i,-\beta)} v.
\end{gather*}
Let $\omega_{\infty} := F_1(\omega_1) \in \omega_{\infty}$. Then $\rho_{+}(s_{ii}(z)) \omega_{\infty} = \omega_{\infty} \begin{cases}
1-za, & i \leq r, \\
1, & i > r.
\end{cases}$

If $r \leq M$ then $\rho_c(s_{ii}(z)) \omega_{\infty} = \omega_{\infty} \begin{cases}
c-zac, & i \leq r, \\
1-zac^2, & i > r.
\end{cases}$ On the other hand for $r > M$ we have $\rho_c(s_{ii}(z)) \omega_{\infty} = \omega_{\infty} \begin{cases}
1-za, & i \leq r, \\
c^{-1}-zac, & i > r.
\end{cases}$
In both $W_{a;c}^{(r)}$ and $L_{r,a}^+$, $\omega_{\infty}$ is the unique (up to scalar multiple) vector annihilated by the $(e_i^+)_{i\in I_0}$ and so by the $s_{ij}(z)$ for $i < j$, because of Proposition~\ref{prop: first properties of inductive system}. A~somewhat surprising observation from Table~\eqref{tab: asym formula} is that $\rho_+(s_{ij}(z)) = 0$ if $r < j < i$.

For the quantum af\/f\/ine superalgebra $U_q(\Gafft)$ in Proposition~\ref{prop: evaluation morphism}, one can def\/ine in the same way the KR modules $W_{k,a}'^{(r)}$ and construct their asymptotic limits $W_{a;c}'^{(r)}, L_{r,a}'^+$.
\begin{defi} \label{def: prefundamental module}The $Y_q(\Glie)$-module $L_{r,a}^-$ is the pullback of the $Y_q(\Glie')$-module $L_{M+N-r,a}'^+$ by the inverse $\CF^{-1}$ of $\CF\colon Y_q(\Glie') \longrightarrow Y_q(\Glie)$ in Proposition~\ref{prop: evaluation morphism}.
\end{defi}As in the positive case, there is a unique (up to scalar multiple) vector $\omega_{\infty} \in L_{r,a}^-$ annihilated by the $s_{ij}(z)$ for $i < j$. $\omega_{\infty}$ is of even parity and $\rho_-(s_{ii}(z)) \omega_{\infty} = \omega_{\infty} \begin{cases}
1, & i \leq r, \\
1-za, & i > r.
\end{cases}$ We shall see that $L_{r,a}^{\pm}$ are irreducible $Y_q(\Glie)$-modules.

\begin{Example}[$\Glie = \mathfrak{gl}(2,1)$ and $r = 1$]
Let $k > 0$. Fix a highest $\ell$-weight vector $u_0$ of $W_{k,a}^{(1)}$. By Theorem~\ref{thm: BKK Schur-Weyl duality}, for $1\leq i \leq k$, there exist unique $u_i, u_i'$ such that $u_0 = (e_1^+)^i u_i$ and $u_i = e_2^+ u_i'$. These together with $u_0$ form a basis of $W_{k,a}^{(1)}$ of weight
\begin{gather*} u_i \in \big(W_{k,a}^{(1)}\big)_{(k-i)\epsilon_1+i\epsilon_2},\qquad u_i' \in \big(W_{k,a}^{(1)}\big)_{(k-i)\epsilon_1+(i-1)\epsilon_2 + \epsilon_3}. \end{gather*}
The structural maps $F_{k,l}$ respect these bases because they commute with $e_1^+$, $e_2^+$.

 Firstly compute the action of $(e_i^{\pm}, K_i)$. By weight grading, $e_1^- u_i = a_i u_{i+1}$ with $a_i \in \BC$ for $0\leq i < k$. From $u_i = e_1^+ u_{i+1}$ and relation~\eqref{rel: relations between e_i and f_j} we obtain
\begin{gather*} a_i u_{i+1} = e_1^- e_1^+ u_{i+1} = e_1^+ e_1^- u_{i+1} - \frac{K_1-K_1^{-1}}{q-q^{-1}} u_{i+1} = \left(a_{i+1}- \frac{q^{k-2i-2}-q^{2i+2-k}}{q-q^{-1}}\right) u_{i+1}. \end{gather*}
By convention $a_k := 0$. This recurrence gives $a_i = \frac{(q^k-q^{2i-k})(q-q^{-2i-1})}{(q-q^{-1})^2}$.

From $e_2^+ e_1^- u_i' = e_1^-e_2^+u_i' = e_1^- u_i = a_i u_{i+1}$ follows also $e_1^- u_i' = a_i u_{i+1}'$. Noting $e_2^+ u_i' = u_i$ and $e_2^- u_i' = 0$ (because of the weight grading), $e_2^- u_i = e_2^- e_2^+ u_i' = \frac{K_2-K_2^{-1}}{q-q^{-1}} u_i' = \frac{q^i - q^{-i}}{q-q^{-1}} u_i'$.
Applying $e_1^+$ to this identity and using $[e_1^+,e_2^-] = 0$, we have $e_1^+ u_i' = \frac{q^{i-1}-q^{1-i}}{q^i - q^{-i}} u_{i-1}'$.

Secondly consider $s_{13}$ and $t_{31}$. In view of relation~\eqref{rel: relations between e_0 and f_i} we have
\begin{gather*}
s_{13} = \frac{s_{11}}{q-q^{-1}}\big(qe_1^+e_2^+ - e_2^+ e_1^+\big),\qquad t_{31} = \big(e_1^-e_2^- - q^{-1} e_2^-e_1^-\big)\frac{t_{11}}{q-q^{-1}}, \\
s_{13}u_i' = \frac{s_{11}}{q-q^{-1}}\left(q-\frac{q^{i-1}-q^{1-i}}{q^i - q^{-i}}\right) u_{i-1} = \frac{q^{k+1}}{q^i - q^{-i}} u_{i-1}, \\
t_{31} u_i = \frac{q^{i-k}}{q-q^{-1}} \left(a_i \frac{q^i - q^{-i}}{q-q^{-1}} - q^{-1}\frac{q^{i+1}-q^{-i-1}}{q-q^{-1}} a_i\right) u_{i+1}' = - \frac{(1-q^{2i-2k})(1-q^{-2i-2})}{(q-q^{-1})^3} u_{i+1}'.
\end{gather*}
Let $E_{ij} \in \End\big(W_{k,a}^{(1)}\big)$ with $i,j \in \{0,1,\dots, k, \overline{1},\dots, \overline{k} \}$ be elementary matrices with respect to the basis $(u_i,u_i')$. In summary, the $U_q(\Gaff)$-module structure on $W_{k,a}^{(1)}$ is given by
\begin{gather*}
\rho^k(s_{11}(z)) = q^k \times \left(\sum_{i=0}^k \big(q^{-i} - zaq^i\big) E_{ii} + \sum_{i=1}^k \big(q^{-i}-zaq^i\big)E_{\overline{i}\overline{i}}\right), \\
\rho^k(s_{22}(z)) = \sum_{i=0}^k \big(q^i - zaq^{2k-i}\big) E_{ii} + \sum_{i=1}^k \big(q^{i-1}-zaq^{2k-i+1}\big) E_{\overline{i}\overline{i}}, \\
\rho^k(s_{33}(z)) = \sum_{i=0}^k \big(1-zaq^{2k}\big) E_{ii} + \sum_{i=1}^k \big(q^{-1}-zaq^{2k+1}\big) E_{\overline{i}\overline{i}},
\\
\rho^k(s_{12}(z)) = q^k \times \left( \sum_{i=0}^{k-1} q^{-i} E_{i,i+1} + \sum_{i=1}^{k-1} \frac{1-q^{-2i}}{q^{i+1}-q^{-i-1}} E_{\overline{i},\overline{i+1}} \right), \\
\rho^k(s_{13}(z)) = q^k \times \sum_{i=1}^k \frac{q}{q^i - q^{-i}} E_{i-1,\overline{i}}, \\
\rho^k(s_{23}(z)) = \sum_{i=1}^k q^i E_{i,\overline{i}},
\\
\rho^k(s_{21}(z)) = za\sum_{i} \frac{(q^{2i}-q^{2k})(q^{i+1}-q^{-i-1})}{(q-q^{-1})^2} (E_{i+1,i} + E_{\overline{i+1},\overline{i}}), \\
\rho^k(s_{31}(z)) = za\sum_{i=0}^{k-1} \frac{(q^{2k}-q^{2i})(1-q^{-2i-2})}{(q-q^{-1})^3} E_{\overline{i+1},i},\\
\rho^k(s_{32}(z)) = za \sum_{i=1}^k \frac{q^{2k-2i}-q^{2k}}{q-q^{-1}} E_{\overline{i},i}.
\end{gather*}
Here for $\rho^k(s_{21}(z))$, the summation $\sum_i$ is understood to be $0\leq i < k$ for $E_{i+1,i}$ and $1\leq i < k$ for $E_{\overline{i+1},\overline{i}}$. Letting $k \rightarrow \infty$ and replacing $(\rho^k, q^k)$ in the above formulas by $(\rho_c, c)$ gives the $U_q(\Gaff)$-module $W_{a;c}^{(1)}$. (Note that $\rho_c(t_{ij}^{(n)}) = - a^{-1}c^{-2} \rho_c\big(s_{ij}^{(1-n)}\big)$.) Dividing $\rho^k(s_{ij}(z))$ by $q^{k}$ whenever $i \leq 1$ and then setting $q^k = 0$, we obtain the $Y_q(\Glie)$-module $L_{1,a}^+$.
\end{Example}
In terms of Young diagrams, the $u_i$, $u_i'$ correspond to the following tableaux (let $k=3$)
\begin{gather*}
u_0 = \young(111),\qquad u_1 = \young(112),\qquad u_2 = \young(122),\qquad u_3 = \young(222), \\
u_1' = \young(113),\qquad u_2' = \young(123), \qquad u_3' = \young(223).
\end{gather*}
\begin{Example}[$\Glie = \mathfrak{gl}(2,1)$ and $r = 2$] \label{example: (2,1,2)}
For $k >0$, let $u_4 \in W_{k,a}^{(2)}$ be such that $e_i^- u_4 = 0$ for $i = 1,2$. It is included in a basis $(u_1,u_2,u_3,u_4)$ of $W_{k,a}^{(2)}$ by Theorem~\ref{thm: BKK Schur-Weyl duality}:
\begin{gather*}
u_3 := e_2^+ u_4,\qquad u_2 := e_1^+ u_3,\qquad u_1 := e_2^+ u_2.
\end{gather*}
From relations \eqref{rel: relations between e_i and f_j}--\eqref{rel: Cartan relations for weight grading}, one deduces the action of the $e_i^{\pm}, K_i$ and then $U_q(\Glie)$. We identify the vector superspace $W_{k,a}^{(2)}$ with $W_{1,a}^{(2)}$ by this basis. The structure maps $F_{k,l}$ are identity maps because they commute with the~$e_i^+$. Let $E_{ij} \in \End\big(W_{k,a}^{(2)}\big)\colon u_k \mapsto u_i \delta_{jk}$. Then
\begin{gather*}
\rho^k(s_{11}(z)) = q^k \times \big((1 - za)(E_{11}+E_{22}) + \big(q^{-1}-zaq\big) (E_{33}+E_{44})\big), \\
\rho^k(s_{22}(z)) = q^k \times \big( (1 - za)(E_{11}+E_{33}) + \big(q^{-1}-zaq\big) (E_{22}+E_{44})\big), \\
\rho^k(s_{33}(z)) = \big(1-zaq^{2k}\big) E_{11} + \big(q^{-1}-zaq^{2k+1}\big)(E_{22}+E_{33}) + \big(q^{-2}-zaq^{2k+2}\big)E_{44},
\\
\rho^k(s_{12}(z)) = q^k \times \big(1-q^{-2}\big) E_{23}, \\
\rho^k(s_{13}(z)) = q^{k} \times \big(q^{-1}-q^{-3}\big) (qE_{24}-E_{13}), \\
\rho^k(s_{23}(z)) = q^{k} \times \big(1-q^{-2}\big) (E_{12} + E_{34}),
\\
\rho^k(s_{21}(z)) = q^k \times za \big(q^2-1\big) E_{32}, \\
\rho^k(s_{31}(z)) = za \big(q^{2k+2}-1\big) E_{42} - za\big(q^{2k+2}-q^2\big) E_{31}, \\
\rho^k(s_{32}(z)) = za \big(q^{2k+1} - q\big)E_{21} + za \big(q^{2k+2}-1\big) E_{43}.
\end{gather*}
The modules $W_{c;a}^{(2)}$, $L_{2,a}^+$ are then obtained as in the previous example.
\end{Example}
Again in terms of tableaux $u_1$, $u_2$, $u_3$, $u_4$ correspond to (let $k=3$)
\begin{gather*} u_1 = \young(111,222),\qquad u_2 = \young(111,223),\qquad u_3 = \young(112,223),\qquad u_4 = \young(113,223). \end{gather*}
\begin{rem} \label{rem: Tsuboi} The $L_{r,a}^{\pm}$ were previously obtained in~\cite{BT,Tsuboi} from the {\it contracted quantum superalgebra} $\dot{\Uc}_q(\Glie)$. It is a superalgebra def\/ined in the same as $U_q(\Glie)$ in Def\/inition~\ref{def: quantum affine superalgebra} except that the~$t_{ii}$ are not required to be invertible. $\ev_a$ in Proposition~\ref{prop: evaluation morphism} degenerates to $\dot{\ev}_a\colon Y_q(\Glie) \longrightarrow \dot{\Uc}_q(\Glie)$. The $L_{r,a}^{\pm}$ are pullbacks of oscillator modules over $\dot{\Uc}_q(\Glie)$ by $\dot{\ev}_a$.
\end{rem}
\begin{lem} \label{lem: finite-dim weight asym}
Let $\beta \in \BQ_+$. The series $\dim \big(W_{k,a}^{(r)}\big)_{k\varpi_r-\beta}$ converges as $k \rightarrow \infty$.
\end{lem}
\begin{proof}
Let $U_q^-(\Glie)$ be the subalgebra of $U_q(\Glie)$ generated by the $(e_i^-)_{i\in I_0}$. Then $U_q(\Glie)_{-\beta}$ is of dimension $c_{\beta} < \infty$. By def\/inition, $U_q^-(\Glie) \longrightarrow W_{k,a}^{(r)}, x \mapsto x \omega_k$ is surjective and sends $U_q(\Glie)_{-\beta}$ to $\big(W_{k,a}^{(r)}\big)_{k\varpi_r-\beta}$. This shows that the series is bounded above by $c_{\beta}$. On the other hand it is increasing by Proposition~\ref{prop: first properties of inductive system}. So it must converge.
\end{proof}

\section[Category $\BGG$ and Baxter's relations]{Category $\boldsymbol{\BGG}$ and Baxter's relations} \label{sec: O and q}
We introduce a category $\BGG$ of $Y_q(\Glie)$-modules including $W_{a;b}^{(r)}$ and $L_{r,a}^{\pm}$ in the spirit of Hernandez--Jimbo~\cite{HJ} and study its Grothendieck ring via $q$-characters of~\cite{FR}.

The quantum af\/f\/ine superalgebra $U_q(\Gaff)$ admits another system of generators, the so-called {\it Drinfeld loop generators}, arising from the Gauss decomposition:
\begin{gather*}
S(z) = \bigg(\sum\limits_{i<j} f_{ji}^+(z) \otimes E_{ji} + 1 \otimes \Id_{\BV}\bigg) \bigg(\sum\limits_l K_l^+(z) \otimes E_{ll} \bigg)\bigg(\sum\limits_{i<j} e_{ij}^+(z) \otimes E_{ij} + 1 \otimes \Id_{\BV}\bigg), \\
T(z) = \bigg(\sum\limits_{i<j} f_{ji}^-(z) \otimes E_{ji} + 1 \otimes \Id_{\BV}\bigg) \bigg(\sum\limits_l K_l^-(z) \otimes E_{ll} \bigg)\bigg(\sum\limits_{i<j} e_{ij}^-(z) \otimes E_{ij} + 1 \otimes \Id_{\BV}\bigg).
\end{gather*}
For example, $K_1^+(z) = s_{11}(z)$ and $K_1^-(z) = t_{11}(z)$. We refer to \cite[Section~3]{Z2} for more details on the relations and on coproduct formulas of these Drinfeld generators. Recall the $d_i$ from equation~\eqref{def: shift}. Def\/ine $\theta_i$ for $i \in I$ inductively by $\theta_1 = 1,\ \theta_{i+1} = q_{i+1}q_i \theta_i$. Def\/ine
\begin{gather*} 
C_i(z) := \prod_{j=1}^i K_j^+(z \theta_j)^{d_j} = \sum_{n\geq 0} C_{i,n}z^n,\qquad K_i^+(z) =: \sum_{n\geq 0} K_{i,n}^+ z^n \in Y_q(\Glie)[[z]].
\end{gather*}
\begin{prop}[{\cite[Theorem 3.5, Proposition 3.6]{Z2}}] \label{prop: JM elements}
Let $k \in I$ and $m,n \geq 0$.
\begin{itemize}\itemsep=0pt
\item[$(1)$] $\Delta (K_{k,m}^+) - \sum\limits_{l=0}^m K_{k,l}^+ \otimes K_{k,m-l}^+ \in \sum\limits_{\alpha \in \BQ_{+} \setminus \{0\}} Y_q(\Glie)_{\alpha} \otimes Y_q(\Glie)_{-\alpha}$.
\item[$(2)$] For all $i,j \in I$ such that $i,j \leq k$, the $s_{ij}^{(n)}$, $t_{ij}^{(n)}$ commute with $C_{i,m}$.
\end{itemize}
\end{prop}
Set $\wCP := (\BC^{\times})^I \times \super,$ and $\lCP := (\BC[[z]]^{\times})^I \times \super$. The multiplicative group structure on~$\BC^{\times}$,~$\BC[[z]]^{\times}$ and the {\it additive} group structure on the ring~$\super$ make~$\wCP$,~$\lCP$ into multiplicative abelian groups. $\wCP$ is naturally a subgroup of~$\lCP$, and $\BC[[z]]^{\times} \longrightarrow \BC^{\times}$, $h(z) \mapsto h(0)$ induces a projection $\varpi\colon \lCP \longrightarrow \wCP$. There is an injective homomorphism of abelian groups{\samepage
\begin{gather*} 
q^?\colon \ \BP \longrightarrow \wCP,\qquad \lambda \mapsto q^{\lambda} := \big(\big(q^{(\epsilon_i,\lambda)}\big)_{i\in I}; |\lambda|\big).
\end{gather*}
We shall view $h(z) \in \BC[[z]]^{\times}$ as the element $((h_i(z) = h(z))_{i\in I};\even)$ in $\lCP$.}

Let $V$ be a $U_q(\Glie)$-module. For $p = ((p_i)_{i\in I};s) \in \wCP$, def\/ine
\begin{gather} \label{def: weight space with parity}
 V_p := \big\{ v \in V_s \, |\, s_{ii} v = p_i v \ \mathrm{for}\ i \in I \big\}.
\end{gather}
If $V_p \neq 0$, then $p$ is called a {\it weight} of $V$, and $V_p$ the weight space of weight~$p$. Let $\wt(V)$ denote the set of weights of $V$. Notice that $U_q(\Glie)_{\alpha} V_p \subseteq V_{q^{\alpha}p}$ for $p \in \wt(V)$ and $\alpha \in \BQ$.

\begin{lem}
For $\lambda, \mu \in \BP$, the weight space $L(\lambda)_{\mu}$ in~\eqref{def: weight space} and $L(\lambda)_{q^{\mu}}$ in~\eqref{def: weight space with parity} coincide.
\end{lem}
\begin{proof}
Clearly $L(\lambda)_{q^{\mu}} \subseteq L(\lambda)_{\mu}$ for all $\mu \in \BP$. Furthermore $L' := \oplus_{\mu \in \BP} L(\lambda)_{q^{\mu}}$ is easily seen to be a sub-$U_q(\Glie)$-module. Since $v_{\lambda} \in L'$ and since $L(\lambda)$ is irreducible, $L(\lambda) = L'$. This implies $L(\lambda)_{q^{\mu}} = L(\lambda)_{\mu}$ for all $\mu \in \BP$.
\end{proof}

Let $V$ be a $Y_q(\Glie)$-module. For $\Bf = ((f_i(z))_{i\in I};s) \in \lCP$ def\/ine
\begin{gather*} V_{\Bf} := \big\{ v \in V_s\, |\, \exists\, d \in \BZ_{>0} \ \mathrm{such\ that}\ (K_{i}^+(z)-f_{i}(z))^d v = 0\ \mathrm{for}\ i \in I \big\}. \end{gather*}
If $V_{\Bf} \neq 0$, then $\Bf$ is called an $\ell$-weight of $V$, and $V_{\Bf}$ the $\ell$-weight space of $\ell$-weight~$\Bf$. Let~$\lwt(V)$ be the set of $\ell$-weights of $V$. As in Section~\ref{sec: KR}, a non-zero $\super$-homogeneous vector $v$ in a $Y_q(\Glie)$-module $V$ is called a {\it highest $\ell$-weight vector} if $s_{jk}^{(n)} v = 0$ for $j < k$ and $s_{ii}(z) v = g_i(z) v$ with $g_i(z) \in \BC[[z]]^{\times}$. By Gauss decomposition we have $K_i^+(z) v = g_i(z) v$, so $v$ is in the $\ell$-weight space of $\ell$-weight $((g_i(z))_{i\in I}; |v|)$. If furthermore $V = Y_q(\Glie) v$, then $V$ is called a {\it highest $\ell$-weight module}, and $((g_i(z))_{i\in I}; |v|)$ the highest $\ell$-weight of~$V$.

By the comments above Def\/inition~\ref{def: prefundamental module}, any one the $Y_q(\Glie)$-modules $W_{k,a}^{(r)}$, $W_{a;c}^{(r)}$, $L_{r,a}^{\pm}$ contains a unique highest $\ell$-weight vector $\omega_{\infty}$, and its $\ell$-weight is denoted by $\Bw_{k,a}^{(r)}$, $\Bw_{a;c}^{(r)}$, $\Bw_{r,a}^{\pm}$ correspon\-ding\-ly. Explicitly, $\Bw_{k,a}^{(r)} = \Bw_{a;q_r^k}^{(r)}$, $\Bw_{a;1}^{(r)} = 1-za$ and
\begin{gather}
\Bw_{r,a}^+ = (\underbrace{1-za,\dots,1-za}_r,\underbrace{1,\dots,1}_{M+N-r};\even), \qquad \Bw_{r,a}^- = (\underbrace{1,\dots,1}_r, \underbrace{1-za,\dots,1-za}_{M+N-r};\even), \label{def: prefund weight} \\
\begin{tabular}{|c|c|c|}
 \hline
 &$r \leq M$ & $r > M$ \\
 \hline
 $\frac{\Bw_{a;c}^{(r)}}{\Bw_{a;1}^{(r)}} $ & $(\underbrace{c,\dots, c}_{r}, \underbrace{g(z),\dots, g(z)}_{M+N-r};\even )$ & $(\underbrace{1,\dots,1}_r, \underbrace{h(z),\dots, h(z)}_{M+N-r}; \even)$ \\
 \hline
 & $g(z) = \frac{1-zac^2}{1-za}$ & $h(z) = \frac{c^{-1}-za c}{1-za}$ \\
 \hline
\end{tabular} \label{def: asym fund weight}
\end{gather}

\begin{defi}Category $\BGG$ is a full subcategory of the category of $Y_q(\Glie)$-modules. An object of $\BGG$ is a $Y_q(\Glie)$-module $V$ subject to the following conditions:
\begin{itemize}\itemsep=0pt
\item[(i)] it has a weight space decomposition $V = \oplus_{p\in \wCP} V_p$ with $\dim V_p < \infty$ for all $p \in \wCP$;
\item[(ii)] there exist $\mu_1,\mu_2,\dots,\mu_d \in \wCP$ such that $\mathrm{wt}(V) \subseteq \cup_{j=1}^d ( q^{\BQ_-}\mu_j)$.
\end{itemize}
\end{defi}

As a f\/irst example, let $p = ((p_i)_{i\in I};s) \in \wCP$ and $h(z) \in \BC[[z]]^{\times}$. There exists a unique $Y_q(\Glie)$-module structure, denoted by $\BC^{h(z)p}$, on the one-dimensional vector superspace of parity $s$ such that $s_{ij}^{(n)} = \delta_{n0}\delta_{ij} h(z) p_i$. Clearly $\lwt\big(\BC^{h(z)p}\big) = \{h(z) p\}$.

By Corollary \ref{cor: generic asym repre} and Lemma \ref{lem: finite-dim weight asym}, the $Y_q(\Glie)$-modules $W_{a;c}^{(r)}, L_{r,a}^{\pm}$ are in category~$\BGG$.

Category $\BGG$ is monoidal (closed under tensor products) and abelian. Any $Y_q(\Glie)$-module $V$ in category $\BGG$ is a direct sum of its $\ell$-weight spaces and has $q$-character~\cite{FR}
\begin{gather} \label{def: q-char}
 \chi_q(V) = \sum_{\Bf \in \mathrm{wt}_{\ell}(V)} \dim (V_{\Bf}) \Bf \in \CEl.
\end{gather}
Here the target $\CEl$ is the set of formal sums $\sum\limits_{\Bf \in \lCP} n_{\Bf} \Bf$ of the symbols $\Bf$ with integer coef\/f\/icients $n_{\Bf}$ such that $\oplus_{\Bf\in \lCP} (\BC^{\varpi(\Bf)})^{\oplus |n_{\Bf}|}$ is in category~$\BGG$. It admits a ring structure: addition is the usual one of formal sums; multiplication is induced by that of~$\lCP$.

Replacing $\mathrm{wt}_{\ell}$ by $\mathrm{wt}$ in equation~\eqref{def: q-char} def\/ines the classical character $\chi(V)$.
\begin{lem} \label{lem: multiplicative structure of q character}
We have $\chi_q(V\otimes W) = \chi_q(V)\chi_q(W)$ for $V,W$ in category~$\BGG$.
\end{lem}
\begin{proof} The proof is the same as in the non-graded case \cite[Remark 2.6]{FR} based on Proposi\-tion~\ref{prop: JM elements}(1) and the partial order on $\wCP$ induced by~$q^{\BQ_+}$.
\end{proof}

\begin{rem} \label{rem: normalized q character} Let $V$ be in category $\BGG$. Suppose there exists $p \in \wCP$ such that the weight space $V_p$ is one-dimensional and $\wt(V) \subseteq q^{\BQ_-}p$. Then $V_p$ is also an $\ell$-weight space of $\ell$-weight~$\Bf$. Def\/ine the normalized character and normalized $q$-character of $V$ by
\begin{gather*}
\widetilde{\chi}(V) := p^{-1} \chi(V),\qquad \widetilde{\chi}_q(V) := \Bf^{-1} \chi_q(V).
\end{gather*}
This is the case when $V$ is any tensor product of the $W_{a;c}^{(r)}$, $L_{r,a}^{\pm}$. If $D$ is a one-dimensional $Y_q(\Glie)$-module, then the normalized $q$-characters of $V$, $V\otimes D$, $D \otimes V$ coincide.
\end{rem}

\begin{rem} \label{rem: q character and tensor product} Let $V,W$ be in category $\BGG$. Let $v \in V$ be a highest $\ell$-weight vector of $\ell$-weight $\Bf \in \lCP$. Then $v \otimes W_{\Bn} \subseteq (V \otimes W)_{\Bn\Bf}$ for $\Bn \in \lwt(W)$. This follows from Proposition~\ref{prop: JM elements}(1); the term $Y_q(\Glie)_{\alpha} v$ for $\alpha \in \BQ_+\setminus \{0\}$ vanishes.
\end{rem}
\begin{lem} \label{lem: character asymptotic modules}
Let $r \in I_0$ and $a,c \in \BC^{\times}$. Then $\widetilde{\chi}_q\big(W_{k,a}^{(r)}\big)$ converges in $\CEl$ as $k \rightarrow \infty$, and
\begin{gather}
\widetilde{\chi}\big(W_{a;c}^{(r)}\big) = \lim_{k\rightarrow \infty} \widetilde{\chi}\big(W_{k,a}^{(r)}\big) = \widetilde{\chi}\big(L_{r,a}^{+}\big) = \chi\big(L_{r,a}^{+}\big), \label{equ: limit normalized characters} \\
\widetilde{\chi}_q\big(W_{a;c}^{(r)}\big) = \lim_{k\rightarrow \infty} \widetilde{\chi}_q\big(W_{k,a}^{(r)}\big) = \widetilde{\chi}_q\big(L_{r,a}^+\big). \label{equ: limit normalized q characters}
\end{gather}
\end{lem}
\begin{proof} We simplify notation $\Bw_{k,a }^{(r)}=: \Bw^k \in \lCP$. Let $\Bn \in \lCP$ with $\Bw^l \Bn \in \lwt(W_{l,a}^{(r)})$. Then $\varpi(\Bn) = q^{-\beta}$ with $\beta \in \BQ_+$. By Remark~\ref{rem: q character and tensor product}, $F_{k,l}\big(W_{l,a}^{(r)}\big)_{\Bw^l\Bn} \subseteq \big(W_{k,a}^{(r)}\big)_{\Bw^k \Bn}$. So the series $\big\{\dim \big(W_{k,a}^{(r)}\big)_{\Bw^k \Bn} \colon k > 0 \big\}$ is increasing and bounded by $c_{\beta}$ in the proof of Lemma~\ref{lem: finite-dim weight asym}. This proves the convergence of the $\widetilde{\chi}_q\big(W_{k,a}^{(r)}\big)$ as $k \rightarrow \infty$.

Next, f\/ix $\Bn = ((n_i(z))_{i\in I}; |\beta|)$ with $\varpi(\Bn) = q^{-\beta}$ and $0 \neq x \in \big(W_{l,a}^{(r)}\big)_{\Bw^l \Bn}$. For $k > l$ write $\Bw^k = \big(\big(w_i^k(z)\big)_{i\in I};\even\big) \in \lCP$. Since the $\ell$-weight space $\big(W_{k,a}^{(r)}\big)_{\Bw^k \Bn}$ is of dimension $\leq c_{\beta}$,
\begin{gather*}\rho^k\big(K_i^+(z) - h(z)\big)^{c_\beta} F_{k,l}(x) = \big(w_i^k(z)n_i(z) - h(z)\big)^{c_{\beta}} F_{k,l}(x) \qquad \mathrm{for}\quad h(z) \in \BC[[z]]. \end{gather*}
Indeed, because of the commutativity of its coef\/f\/icients, $\rho^k(K_i^+(z))$ restricted to $\big(W_{k,a}^{(r)}\big)_{\Bw^k \Bn}$ can be made into an upper triangular matrix with uniform diagonals~$w_i^k(z)n_i(z)$.

Let us take $h(z)$ to be $n_i(z)$ times the $i$-th component $w_i(z)$ of $\Bw_{a;c}^{(r)}$:
\begin{gather*} \rho^k\big(K_i^+(z) - w_i(z) n_i(z)\big)^{c_{\beta}} F_{k,l}(x) = \big(w_i^k(z)-w_i(z)\big)^{c_{\beta}}n_i(z)^{c_{\beta}} F_{k,l}(x). \end{gather*}
At the right-hand side the factor before $F_{k,l}(x)$, when expanded at $z = 0$, has as coef\/f\/icients Laurent polynomials in $q_r^k$. To compute $\rho_c(K_i^+(z) - w_i(z) n_i(z))^{c_{\beta}} F_{l}(x)$ in $W_{a;c}^{(r)}$, it suf\/f\/ices to evaluate these polynomials at $q_r^k = c$ by Remark \ref{rem: asym computation}. But by def\/inition the factor $(w_i^k(z)-w_i(z))|_{q_r^k = c} = 0$. So $\rho_c(K_i^+(z) - w_i(z) n_i(z))^{c_{\beta}} F_{l}(x) = 0$, meaning that $F_l(x) \in W_{a;c}^{(r)}$ is in the $\ell$-weight space of $\ell$-weight $ \Bw_{a;c}^{(r)} \Bn$. This proves the f\/irst equality of equation~\eqref{equ: limit normalized q characters}. The second equality for $L_{r,a}^+$ can be proved in the same way using $\tilde{\rho}^k$.
\end{proof}

Let $\CR$ be the subset of $\lCP$ consisting of the $\Bf = ((f_i(z))_{i\in I};s)$ such that $\frac{f_i(z)}{f_{i+1}(z)}$ is the Taylor expansion at $z=0$ of a rational function for $i\in I_0$.
\begin{lem} \label{lem: simple modules in category O}
Let $\Bf = ((f_i(z))_{i\in I};s) \in \CR$.
\begin{itemize}\itemsep=0pt
\item[$(1)$] There exists a unique (up to isomorphism) irreducible $Y_q(\Glie)$-module in category $\BGG$ of highest $\ell$-weight $\Bf$. Let $V(\Bf)$ be the $Y_q(\Glie)$-module thus obtained.
\item[$(2)$] $V(\Bf)$ can be extended to a $U_q(\Gaff)$-module if and only if for all $i \in I_0$, as a rational function $\frac{f_i(z)}{f_{i+1}(z)}$ is a product of the $\frac{1-za}{c^{-1}-zac}$ with $a,c \in \BC^{\times}$.
\item[$(3)$] All irreducible $Y_q(\Glie)$-modules in category $\BGG$ are of the form $V(\Bf)$ with $\Bf \in \CR$.
\end{itemize}
\end{lem}
\begin{proof}
 (1) and suf\/f\/iciency of (2). In view of equations~\eqref{def: prefund weight}--\eqref{def: asym fund weight}, such an $\Bf$ can be written as $h(z) p \Bn$ where $h(z) \in \BC[[z]]^{\times}$, $p \in \wCP$ (resp. $p$ is of the form $(1,\dots,1;s) \in \wCP$) and $\Bn$ is a~product of the $\Bw_{r,a}^{\pm}$ (resp.\ the $\Bw_{a;c}^{(r)}$). So $V(\Bf)$ can be realized as a sub-quotient of the tensor product of $V(h(z)p)$, which is a one-dimensional $Y_q(\Glie)$-module (resp.\ $U_q(\Gaff)$-module), with the corresponding tensor product of the~$L_{r,a}^{\pm}$ (resp.\ the $W_{a;c}^{(r)}$). This shows that~$V(\Bf)$ is in category~$\BGG$ (resp.\ a~$U_q(\Gaff)$-module).

Necessity of (2). One considers the action of $K_{i+1}^{\pm}(z)K_{i}^{\pm}(z)^{-1}$, based on the Drinfeld relations involving $[X_i^+,X_i^-]$ in \cite[Theorem~3.5]{Z2} and the assumption $\dim V(\Bf)_{\varpi(\Bf) q^{-\alpha_i}} < \infty$. As in \cite[Proposition~6.1, Lemma~4.12]{Z}, $g_i(z) := \frac{f_i(z)}{f_{i+1}(z)}$ is regular at $z = 0, \infty$ and $g_i(0)g_i(\infty) = 1$. Necessarily $g_i(z)$ is a product of the $\frac{1-za}{c^{-1}-zac}$ with $a,c \in \BC^{\times}$.

Similar arguments can be used to prove (3); since $K_i^-(z) \notin Y_q(\Glie)[[z^{-1}]]$, one loses the regularity of $g_i(z)$ at $z = \infty$. See also \cite[Lemma~3.9]{HJ}.
\end{proof}

The abelian category $\BGG$ contains modules with inf\/inite Jordan--H\"older series, so we need its {\it completed} Grothendieck group $K_0(\BGG)$: elements are formal sums $\sum_{\Bf \in \CR} n_{\Bf} [V(\Bf)]$ of the symbols~$[V(\Bf)]$ with integer coef\/f\/icients $n_{\Bf}$ such that $\oplus_{\Bf \in \CR} V(\Bf)^{\oplus |n_{\Bf}|}$ is in category $\BGG$; addition is the usual one of formal sums. As in the case of Kac--Moody algebras \cite[Section~9.3]{Ka}, for $\Bf \in \CR$, the multiplicity $m_{\Bf,X} \in \BZ_{\geq 0}$ of $V(\Bf)$ in any object $X$ of category $\BGG$ is well-def\/ined, and $[X] := \sum\limits_{\Bf \in \CR} m_{\Bf,X} [V(\Bf)] \in K_0(\BGG)$. There is no ambiguity for $X = V(\Bf)$ as $m_{\Bn,V(\Bf)} = \delta_{\Bn,\Bf}$ for $\Bn, \Bf \in \CR$. Make $K_0(\BGG)$ into a ring with multiplication induced by $[X][Y] := [X \otimes Y]$ for $X, Y$ in category $\BGG$. Since $\chi_q$ respects exact sequences and tensor products, the assignment $[X] \mapsto \chi_q(X)$, for $X$ in category $\BGG$, extends uniquely to a ring homomorphism $\chi_q\colon K_0(\BGG) \longrightarrow \CEl$, called the $q$-character map.
\begin{cor} \label{cor: injectivity q-char}
The $q$-character map $\chi_q$ is injective.
\end{cor}
\begin{proof}We need to show that $\chi_q(V(\Bf))$ distinguishes $\Bf$. Indeed, from the proof of Lem\-ma~\ref{lem: simple modules in category O}(1), we deduce that $\chi_q (V(\Bf))$ is $\Bf$ plus terms of the form $\Bn\Bf$ where the $\Bn \in \lCP$ satisfy $\varpi(\Bn) \in q^{\BQ_-}$. So $\Bf$ appears in $\chi_q(V(\Bf))$ as a leading term.
\end{proof}

\begin{prop}Let $\Bf = ((f_i(z))_{i\in I};s) \in \CR$. Then $V(\Bf)$ is finite-dimensional if and only if for all $i \in I_0 \setminus \{M\}$ there exist $P_i(z) \in 1 + z\BC[z]$ and $a_i \in \BC^{\times}$ such that $\frac{f_i(z)}{f_{i+1}(z)} = a_i \frac{P_i(z q_i^{-1})}{P_i(z q_i)}$.
\end{prop}
\begin{proof}(Sketch\footnote{This result is not needed in the following. We include it here for completeness.}) Suf\/f\/iciency: such $\Bf$ can be written as $h(z) p \Bn$ where $h(z) \in \BC[[z]]^{\times}$, $p \in \wCP$ and $\Bn$ is a product of the $\Bw_{a;q_r}^{(r)}$, $\Bw_{M,a}^{\pm}$ with $a \in \BC^{\times}$ and $r \in I_0$. So~$V(\Bf)$ is a sub-quotient of the tensor product $T$ the $V(h(z)p)$, with the $V\big(\Bw_{a;q_r}^{(r)}\big) = W_{1,a}^{(r)}$ and $L_{M,a}^{\pm}$. By Lemma~\ref{lem: character asymptotic modules} and Theorem~\ref{thm: BKK Schur-Weyl duality}, $\dim L_{M,a}^{\pm} = 2^{MN}$. So~$T$ and~$V(\Bf)$ are f\/inite-dimensional.

Necessity: let $i \in I_0$. One restricts to the subalgebra $Y_i$ of $Y_q(\Glie)$ generated by the $s_{jk}^{(n)}$, $s_{jj}^{-1}$ with $j,k \in \{i,i+1\}$ and $n \geq 0$. It is a quotient algebra of $Y_q(\mathfrak{gl}_2)$. The polynomiality can then be proved along the line of \cite[Remark~3.11]{FH}.
\end{proof}

The $L_{r,a}^{\pm}$ in \cite{HJ} are always inf\/inite-dimensional. In our case, by Theorem~\ref{thm: BKK Schur-Weyl duality} and equation~\eqref{equ: limit normalized characters}, $L_{M,a}^+$ is of dimension $2^{MN}$; see also Example~\ref{example: (2,1,2)}. We refer to \cite[Section~4]{Z2} for a~detailed discussion of f\/inite-dimensional irreducible $Y_q(\mathfrak{gl}(1,1))$-modules.

Let $\CR_U$ be the subset of $\CR$ formed of $\Bf$ satisfying the condition in Lemma \ref{lem: simple modules in category O}(2). The following elements of $\CR_U$ will be used: let $r \in I_0$ and $a \in \BC^{\times}$,
\begin{gather} \label{def: generalized simple roots}
A_{r,a} = \left(\underbrace{1,\dots,1}_{r-1}, \frac{q_r - za\theta_r^{-1}}{1 - za\theta_r^{-1}q_r}, \frac{1 - za\theta_{r+1}^{-1}q_{r+1}}{q_{r+1} - za\theta_{r+1}^{-1}}, \underbrace{1,\dots,1}_{M+N-r-1}; |\alpha_r|\right).
\end{gather}
By def\/inition $\varpi(A_{r,a}) = q^{\alpha_r}$. The $A_{r,a}$ generate a free abelian subgroup $\lCQ$ of $\CR_U$.

\begin{thm} \label{thm: Baxter}
Let $S$ be an irreducible $U_q(\Gaff)$-module in category $\BGG$.
\begin{itemize}\itemsep=0pt
\item[$(1)$] We have $\lwt(S) \subseteq \lCQ \Bf$ where $\Bf \in \CR_U$ is the highest $\ell$-weight of $S$.
\item[$(2)$] If $S$ is finite-dimensional, then in a fraction ring of $K_0(\BGG)$,
\begin{gather} \label{rel: TQ}
 [S] = \sum_{j=1}^{\dim S} [D_j] \Bm_j,
\end{gather}
where for each $j$, $D_j$ is a one-dimensional $U_q(\Gaff)$-module in category $\BGG$, and $\Bm_j$ is a product of the $\frac{[W_{a;b}^{(r)}]}{[W_{a;c}^{(r)}]}$ with $r \in I_0$ and $a,b,c \in \BC^{\times}$.
\end{itemize}
\end{thm}
(2) can be though of as generalized Baxter's relations in the sense of Frenkel--Hernandez \cite[Theorem 4.8]{FH}. In equation~\eqref{rel: TQ}, only $U_q(\Gaff)$-modules are involved. Proof of (1) and concrete examples are postponed to Section~\ref{sec: GT bases}.

\begin{proof}[Proof of Theorem \ref{thm: Baxter}(2) assuming (1)] Any $\Bn \in \CR_U$ can be written as $h(z) p \Bn'$ where $h(z) \in \BC[[z]]^{\times}$, $p = (1,\dots,1;s) \in \wCP$ and $\Bn'$ is a product of the $\frac{\Bw_{a;b}^{(r)}}{\Bw_{a;c}^{(r)}} = \frac{\chi_q(W_{a;b}^{(r)})}{\chi_q(W_{a;c}^{(r)})}$; the last identity follows from equation~\eqref{equ: limit normalized q characters}. By (1), $\lwt(S) \subset \CR_U$. This establishes the $q$-character version of equation~\eqref{rel: TQ} with isomorphism classes replaced by $q$-characters ($D_j = V(h(z)p)$). We conclude by the injectivity of the $q$-character map.
\end{proof}

\begin{cor} \label{cor: simplicity prefund}
The $Y_q(\Glie)$-modules $L_{r,a}^{\pm}$ for $r \in I_0$ and $a \in \BC^{\times}$ are irreducible.
\end{cor}
\begin{proof}We show $L_{r,a}^+$ is irreducible, imitating the proof of \cite[Theorem~6.1]{HJ}. By Def\/inition \ref{def: prefundamental module} this implies the irreducibility of $L_{r,a}^-$. Let $S_{r,a}^+$ be the sub-module of $L_{r,a}^+$ generated by its unique highest $\ell$-weight vector. Then $S_{r,a}^+ \cong V(\Bw_{r,a}^+)$ is irreducible.

We prove that $\widetilde{\chi}_q (L_{r,a}^+) = \widetilde{\chi}_q(S_{r,a}^+)$. In view of equation~\eqref{equ: limit normalized q characters}, this means that: if $l > 0$ and $\Bw_{l,a}^{(r)} \Bf \in \lwt\big(W_{l,a}^{(r)}\big)$, then $\dim (S_{r,a}^+)_{\Bw_{r,a}^+\Bf} \geq \dim \big(W_{l,a}^{(r)}\big)_{\Bw_{l,a}^{(r)} \Bf}$. We shall indicate the dependence of an element $\Bm = \Bm(z) \in \lCP$ on $z$ whenever necessary. For $\Bf = 1$ this is clear as they are both one-dimensional. Let $\Bf \neq 1$. By Theorem~\ref{thm: Baxter}(1) and equation~\eqref{equ: limit normalized q characters}, $\Bf \in \lCQ$ and $\varpi(\Bf) \in q^{\BQ_-}$. The following set $X$ is f\/inite and non-empty:
\begin{gather*} X := \big\{ \Bn^{-1}\Bf \in \lCQ \, |\, \Bf \neq \Bn \in \lCQ,\ \Bw_{r,a}^+\Bn \in \lwt(S_{r,a}^+),\ \varpi\big( \Bn^{-1}\Bf\big) \in q^{\BQ_-} \big\} \subset \lCQ \setminus \{1\}.\end{gather*}
For $\Bm \in X$, the set $\{ \Bm(zq_r^{-2k})\, |\, k > l\} $ is inf\/inite but its intersection with $\lwt(L_{r,a}^-) (\Bw_{r,a}^-)^{-1}$ is f\/inite, because $\varpi(\Bm(zq_r^{-2k})) = \varpi(\Bm)$ is f\/ixed and weight spaces of $L_{r,a}^-$ are f\/inite-dimensional. Choose $k>l$ large enough such that $ \Bm(zq_r^{-2k})\Bw_{r,a}^- \notin \lwt(L_{r,a}^-)$ for all $\Bm \in X$. By Lemma~\ref{lem: character asymptotic modules}, $\Bw_{k,a}^{(r)}\Bf \in \lwt(W_{k,a}^{(r)})$. The $Y_q(\Glie)$-module $S_{r,a}^+ \otimes L_{r,aq_r^{2k}}^-$ has an irreducible sub-quotient isomorphic to $W_{k,a}^{(r)}$ up to tensor product by one-dimensional $Y_q(\Glie)$-modules. It follows from Lemma~\ref{lem: multiplicative structure of q character} and Remark~\ref{rem: normalized q character} that $\Bf = \Bf_k^+ \Bf_k^-$ where $\Bw_{r,a}^+ \Bf_k^+ \in \lwt(S_{r,a}^+)$ and $\Bw_{r,aq_r^{2k}}^- \Bf_{k}^- \in \lwt(L_{r,aq_r^{2k}}^-)$, which implies $\Bf_k^-(zq_r^{-2k}) \Bw_{r,a}^- \!\in\! \lwt(L_{r,a}^-)$ based on the pullback of~$\Phi_{q_r^{2k}}$ in Proposition \ref{prop: evaluation morphism}. So \smash{$( \Bf_k^+)^{-1} \Bf = \Bf_k^-\!\notin\! X$}. By Theorem~\ref{thm: Baxter}(1) and equation~\eqref{equ: limit normalized q characters}, $\Bf_k^+ \in \lCQ$ and $\varpi(\Bf_k^-) \in q^{\BQ_-}$, forcing $\Bf_k^- =1$. So any such factorization $\Bf = \Bf_k^+ \Bf_k^-$ is trivial: $\Bf_k^- = 1$. This proves $\dim \big(W_{l,a}^{(r)}\big)_{\Bw_{l,a}^{(r)} \Bf} \leq \dim \big(W_{k,a}^{(r)}\big)_{\Bw_{k,a}^{(r)} \Bf} \leq \dim (S_{r,a}^+)_{\Bw_{r,a}^+\Bf}$.
\end{proof}

\section[Rationality of $\ell$-weights of f\/inite-dimensional modules]{Rationality of $\boldsymbol{\ell}$-weights of f\/inite-dimensional modules} \label{sec: GT bases}
In this section, we study in more detail the $\ell$-weights of f\/inite-dimensional $U_q(\Gaff)$-modules in category $\BGG$ and prove Theorem~\ref{thm: Baxter}(1).

Recall the $U_q(\Glie)$-module $\BV$ def\/ined by \eqref{def: vector repr}. Following \cite[Example~2]{Z3}, def\/ine the $U_q(\Glie)$-module $\BW$ to be the pullback of $\BV$ by $\Psi\colon U_q(\Glie) \longrightarrow U_q(\Glie)$ in Proposition~\ref{prop: evaluation morphism}. For $i \in I$, set $u_i := \Psi^* v_i \in \BW$; it is of weight $q^{-\epsilon_i}$.

For $a \in \BC^{\times}$, def\/ine the $U_q(\Gaff)$-modules $\BV(a)$ and $\BW(a)$ to be the pullbacks of $\BV$ and $\BW$ respectively by $ \ev_a\colon U_q(\Gaff) \longrightarrow U_q(\Glie)$. Naturally as $Y_q(\Glie)$-modules, $\BV(a)$ and~$\BW(a)$ are in category~$\BGG$, and their weight spaces are always one-dimensional.
\begin{defi} \label{defi: fundamental weights}
For $a \in \BC^{\times}$ and $i \in I$, let $\CX_{i,a} \in \lwt(\BV(a))$ and $\CY_{i,a} \in \lwt(\BW(a))$ be such that $\varpi(\CX_{i,a}) = q^{\epsilon_i}$ and $\varpi(\CY_{i,a}) = q^{-\epsilon_i}$.
\end{defi}
Let us compute explicitly the $\CX_{i,a}$, $\CY_{i,a}$, following an idea of \cite[Lemma~4.7]{FM}. Fix $k,i \in I$. By Def\/inition \ref{defi: fundamental weights}, $v_i$, $u_i$ are common eigenvectors of~$C_k(z)$, whose eigenvalues are denoted by $g_i^k(z), h_i^k(z) \in \BC[[z]]^{\times}$ respectively. An important observation from the Gauss decomposition in Section~\ref{sec: O and q} is that: if a vector $x$, either in $\BV(a)$ or $\BW(a)$ is annihilated by the $s_{jl}^{(n)}$, $t_{jl}^{(n)}$ for $1\leq j<l \leq k$, then $K_j^+(z) x = s_{jj}(z) x$ for $1\leq j \leq k$. For example, $v_1$, $u_k$, $v_l$, $u_l$ for $l > k$ are such vectors according to the weight gradings on $\BV$, $\BW$. Therefore for $i > k$,
\begin{gather*} g_i^k(z) = \prod_{j=1}^k (1-za\theta_j)^{d_j} = h_i^k(z). \end{gather*}
Suppose $i \leq k$. Observe from \eqref{def: vector repr} that $v_i$ is proportional to $t_{i1} v_1$. Similarly $u_i$ is proportional to $t_{ki} u_k$; see \cite[Example 2]{Z3}. Since $t_{i1}$, $t_{ki}$ commute with $C_k(z)$ by Proposition~\ref{prop: JM elements}, we have $g_i^k(z) = g_1^k(z)$ and $h_i^k(z) = h_k^k(z)$. We apply the above observation to the vectors $v_1$, $u_k$ to compute $g_1^k(z)$ and $h_k^k(z)$, and obtain
\begin{gather*} g_i^k(z) = \big(q_1-zaq_1^{-1}\big)^{d_1} \prod_{j=2}^k (1-za\theta_j)^{d_j},\qquad h_i^k(z) = \big(q_k^{-1}-za\theta_kq_k\big)^{d_k} \prod_{j=1}^{k-1} (1-za\theta_j)^{d_j}. \end{gather*}

It follows that: setting $h(z) = (1-zaq^2)(1-zaq^{-2})(1-za)^{-2}$,
\begin{gather*}
\CX_{i,a} = (1-za) \times \left( \underbrace{1,\dots,1}_{i-1}, \left(\frac{q-za\theta_i^{-1}q^{-1}}{1-za\theta_i^{-1}}\right)^{d_i}, \underbrace{1,\dots,1}_{M+N-i}; |\epsilon_i| \right), \\
\CY_{i,a} = (1-za) \times \left(\underbrace{1,\dots,1}_{i-1}, \frac{q_i^{-1}-zaq_i}{1-za}, \underbrace{h(z) ,\dots,h(z)}_{M+N-i}; |\epsilon_i| \right).
\end{gather*}
We used $\theta_{i+1} = \theta_i q_iq_{i+1}$ for $i \in I_0$. Combining with equation~\eqref{def: generalized simple roots}, we obtain
\begin{gather} \label{equ: A, X, Y}
A_{r,a} = \CX_{r,aq} \CX_{r+1,aq}^{-1} = \CY_{r,a\theta_r^{-1}q_r^{-1}}^{-1} \CY_{r+1,a\theta_r^{-1}q_r^{-1}}.
\end{gather}

KR modules $W_{1,a}^{(r)}$ can be constructed from $\BV, \BW$ by a fusion procedure.
\begin{lem}[\cite{Z3}] \label{lem: fusion}
Let $1\leq s \leq M$ and $1\leq t < N$. Let $a \in \BC^{\times}$. Then as $U_q(\Gaff)$-modules,
\begin{gather*}
W_{1,a}^{(s)} \simeq U_q(\Gaff)v_{M+N}^{\otimes s} \subseteq \bigotimes_{1\leq i\leq s}^{\leftarrow} \BV\big(aq^{2i}\big),\qquad W_{1,a}^{(M+N-t)} \simeq U_q(\Gaff)u_1^{\otimes t} \subseteq \bigotimes_{1\leq j \leq t}^{\leftarrow} \BW\big(aq^{2t-2j}\big).
\end{gather*}
\end{lem}
\begin{proof}
Comparing Def\/inition \ref{def: KR} with \cite[Lemma 5, Def\/inition 2]{Z3} we have
\begin{gather*}W_{1,a}^{(s)} = \ev_{aq^2}^* L(\varpi_s) \simeq V_{s,aq^{2s+2}}^+ := U_q(\Gaff) v_{M+N}^{\otimes s} \subset \bigotimes_{1\leq i\leq s}^{\leftarrow} \BV\big(aq^{2i}\big). \end{gather*}
Here we borrow the notation $V_{s,aq^{2s+2}}^+$ from {\it loc.\ cit.} The case of $t$ is parallel.
\end{proof}

Note that in equation~\eqref{equ: A, X, Y}, $\theta_r^{-1}q_r^{-1} \in q^{2\BZ + 1}$. Let us def\/ine $\lCQ_a$ to be the subgroup of $\lCQ$ generated by the $A_{r,b}$ with $b \in aq^{2\BZ}$ and $r \in I_0$.
\begin{cor} \label{cor: rationality} For $r \in I_0, a,c \in \BC^{\times}$ and $k \in \BZ_{>0}$ we have
\begin{gather*} \big(\Bw_{k,a}^{(r)}\big)^{-1} \lwt\big(W_{k,a}^{(r)}\big) \subseteq \big(\Bw_{r,a}^+\big)^{-1} \lwt(L_{r,a}^+) = \big(\Bw_{a;c}^{(r)}\big)^{-1} \lwt\big(W_{a;c}^{(r)}\big) \subseteq \lCQ_{aq}. \end{gather*}
\end{cor}
\begin{proof}By Lemmas \ref{lem: KR fund modules} and~\ref{lem: fusion}, there exists a one-dimensional $U_q(\Gaff)$-module $D$ such that $D \otimes W_{k,a}^{(r)}$ is a sub-quotient of a tensor product $T$ of the $X(b)$ where $b \in aq^{2\BZ}$ and: if $r \leq M$ then $X = \BV$; if $r > M$ then $X = \BW$. For any such $X(b)$, let $\Bf$ be the highest $\ell$-weight of $X(b)$. Then by equation~\eqref{equ: A, X, Y}: $\Bf^{-1} \lwt(X(b)) \subset \lCQ_{aq}$. Since $\chi_q$ respects tensor products, we have $\Bn^{-1}\Bm \in \lCQ_{aq}$ for all $\Bn,\Bm \in \lwt(T)$. Taking normalized $q$-characters in Remark~\ref{rem: normalized q character}, $\big(\Bw_{k,a}^{(r)}\big)^{-1} \lwt\big(W_{k,a}^{(r)}\big) \subseteq \lCQ_{aq}$. The rest comes from equation~\eqref{equ: limit normalized q characters}.
\end{proof}

\begin{proof}[Proof of Theorem \ref{thm: Baxter}(1)] As in the proof of suf\/f\/iciency of Lemma~\ref{lem: simple modules in category O}(2), the irreducible module $V(\Bf)$ for $\Bf \in \CR_U$ can be realized as a sub-quotient of a tensor product $T$ of a one-dimensional $U_q(\Gaff)$-module $D$ with the $W_{c;a}^{(r)}$. Comparing normalized $q$-characters, we obtain (1) as a consequence of Corollary~\ref{cor: rationality}.
\end{proof}

Another consequence is the irreducibility of $W_{a;c}^{(r)}$ for generic $c \in \BC^{\times}$.
\begin{cor}The $U_q(\Gaff)$-module $W_{a;c}^{(r)}$ is irreducible if $c^2 \notin q^{2\BZ}$.
\end{cor}
\begin{proof}Let $S_{a;c}^{(r)}$ be the sub-module of $W_{a;c}^{(r)}$ generated by its highest $\ell$-weight vector; it is irreducible and isomorphic to $V\big(\Bw_{a;c}^{(r)}\big)$. Let $k > 0$ and $\Bw_{k,a}^{(r)}\Bf \in \lwt\big(W_{k,a}^{(r)}\big)$. By equation~\eqref{def: asym fund weight}
\begin{gather*} \big(1-zac^2\big) \times \Bw_{k,a}^{(r)} = \Bw_{a;c}^{(r)} \Bw_{ac^2;c^{-1}q_r^k}^{(r)} \in \lCP. \end{gather*}
This implies that $V(1-zc^2) \otimes W_{k,a}^{(r)}$ is an irreducible sub-quotient of $S_{a;c}^{(r)} \otimes W_{ac^2;c^{-1}q_r^k}^{(r)}$. Taking normalized $q$-characters and noting that $V(1-zac^2)$ is one-dimensional, we have $\Bf = \Bf_+ \Bf_-$ where $\Bw_{a;c}^{(r)} \Bf_+ \in \lwt\big(W_{a;c}^{(r)}\big)$ and $\Bw_{ac^2;c^{-1}q_r^k}^{(r)} \Bf_- \in \lwt\big(W_{ac^2;c^{-1}q_r^k}^{(r)}\big)$. By Corollary~\ref{cor: rationality}, we have $\Bf, \Bf_+ \in \lCQ_{aq}$ and $\Bf_- \in \lCQ_{ac^2q}$. The assumption on $c$ implies that $\lCQ_{aq}\cap \lCQ_{ac^2q}=\{1\}$, which forces $\Bf_- = 1$. So any such factorization $\Bf = \Bf_+\Bf_-$ is trivial: $\Bf_- = 1$. This shows that $\dim \big(W_{k,a}^{(r)}\big)_{\Bw_{k,a}^{(r)}\Bf} \leq \dim \big(S_{a;c}^{(r)}\big)_{\Bw_{a;c}^{(r)} \Bf}$ for all $k$, and therefore $S_{a;c}^{(r)} = W_{a;c}^{(r)}$.
\end{proof}

\begin{Example}[Baxter's relations for $\mathfrak{gl}(2,1)$] We derive equation~\eqref{rel: TQ} for the $U_q(\Gaff)$-modules $\BV(a)$ and $\BW(a)$, along the line of the proof of Theorem~\ref{thm: Baxter}(2):
\begin{gather*}
\frac{\Bw_{a;c}^{(1)}}{\Bw_{a;1}^{(1)}} = \left(c, \frac{1-zac^2}{1-za}, \frac{1-zac^2}{1-za};\even \right),\qquad \frac{\Bw_{a;c}^{(2)}}{\Bw_{a;1}^{(2)}} = \left(c,c, \frac{1-zac^2}{1-za};\even \right), \\
\CX_{1,a} = (1-za) \times \left( \frac{q-zaq^{-1}}{1-za},1,1; \even \right) = \big(1-zaq^{-2}\big) \frac{\Bw_{aq^{-2};q}^{(1)}}{\Bw_{aq^{-2};1}^{(1)}}, \\
\CX_{2,a} = (1-za) \times \left(1, \frac{q-zaq^{-3}}{1-zaq^{-2}}, 1;\even \right) = (1-za) \frac{\Bw_{aq^{-2};q^{-1}}^{(1)}}{\Bw_{aq^{-2};1}^{(1)}} \frac{\Bw_{aq^{-4};q}^{(2)}}{\Bw_{aq^{-4};1}^{(2)}}, \\
\CX_{3,a} = (1-za) \left(1,1, \frac{1-zaq^{-2}}{q-zaq^{-3}} ; \odd \right) = (1-za) \frac{\Bw_{aq^{-4};q}^{(2)}}{\Bw_{aq^{-4};1}^{(2)}} \big(q^{-1},q^{-1},q^{-1};\odd\big), \\
\CY_{1,a} = \left(q^{-1}-zaq, \frac{(1-zaq^2)(1-zaq^{-2})}{1-za}, \frac{(1-zaq^2)(1-zaq^{-2})}{1-za}; \even \right)\\
 \hphantom{\CY_{1,a}}{} = \big(1-zaq^2\big)\frac{\Bw_{a;q^{-1}}^{(1)}}{\Bw_{a;1}^{(1)}}, \\
\CY_{2,a} = \left(1-za,q^{-1}-zaq, \frac{(1-zaq^2)(1-zaq^{-2})}{1-za};\even\right) = (1-za) \frac{\Bw_{a;q}^{(1)}}{\Bw_{a;1}^{(1)}} \frac{\Bw_{a;q^{-1}}^{(2)}}{\Bw_{a;1}^{(2)}}, \\
\CY_{3,a} = (1-za) \times \left(1,1,\frac{q-zaq^{-1}}{1-za};\odd \right) = (1-za) \frac{\Bw_{a;q^{-1}}^{(2)}}{\Bw_{a;1}^{(2)}} (q,q,q;\odd).
\end{gather*}
The irreducible $U_q(\Gaff)$-module of highest $\ell$-weight $(g(z),g(z),g(z);s)$ being one-dimensional for $g(z) \in \BC[[z]]^{\times}$ and $s \in \super$, its isomorphism class is denoted by $[g(z)s]$. Then
\begin{gather*}
\big[\BV\big(aq^2\big)\big] = [1-za] \frac{[W_{a;q}^{(1)}]}{[W_{a;1}^{(1)}]} + [1-zaq^2]\frac{[W_{a;q^{-1}}^{(1)}]}{[W_{a;1}^{(1)}]} \frac{[W_{aq^{-2};q}^{(2)}]}{[W_{aq^{-2};1}^{(2)}]} + [(q^{-1}-zaq)\odd] \frac{[W_{aq^{-2};q}^{(2)}]}{[W_{aq^{-2};1}^{(2)}]}, \\
[\BW(a)] = [1-zaq^2]\frac{[W_{a;q^{-1}}^{(1)}]}{[W_{a;1}^{(1)}]} + [1-za] \frac{[W_{a;q}^{(1)}]}{[W_{a;1}^{(1)}]} \frac{[W_{a;q^{-1}}^{(2)}]}{[W_{a;1}^{(2)}]} + [(q-zaq)\odd] \frac{[W_{a;q^{-1}}^{(2)}]}{[W_{a;1}^{(2)}]}.
\end{gather*}
\end{Example}

To derive Baxter's relations, one needs to compute its $q$-character of f\/inite-dimensional $U_q(\Gaff)$-modules. In a previous version of this paper (arXiv:1410.0837v2), the author was able to prove a tableau-sum formula for $\chi_q(\ev_a^*L(\lambda))$, with $L(\lambda)$ being an irreducible submodule of a tensor power of~$\BV$ (polynomial modules), based on an idea of~\cite{FM} relating $\ell$-weights to Gelfand--Tsetlin basis~\cite{PSV}. Such formula appeared earlier in the context of transfer matrices of quantum integrable systems attached to~$U_q(\Gaff)$~\cite{Tsuboi1}. In a recent preprint~\cite{Z4}, the tableau-sum $q$-character formula has been extended to
\begin{gather*}\ev_a^* (L^*(\lambda)),\qquad (\ev_a')^*L(\lambda),\qquad (\ev_a')^*(L^*(\lambda)),\end{gather*}
where $\ev_a'$ is a second evaluation map arising from involution and $L^*(\lambda)$ is the dual module. To avoid redundancy, in the present paper we do not present the $q$-character formula.

\appendix

\section{Non-twisted quantum loop algebras}
We apply the asymptotic construction ($c \in \BC^{\times}$) of Section~\ref{sec: idea} to the inductive system in \cite{HJ} of KR modules over an arbitrary non-twisted quantum loop algebra.

The main step is to establish the asymptotic property as Lemma \ref{lem: asym property KR}. In Section~\ref{sec: generic asymptotic representations}, we used the evaluation maps to reduce to the f\/inite-type quantum supergroup (Proposition \ref{prop: first properties of inductive system}). What is essential is the representation theory of $U_q(\mathfrak{sl}_2)$ and $U_q(\mathfrak{sl}(1,1))$. The evaluation maps, which do not exist for quantum loop algebras out of type~A, do not play big r\^{o}le.

We use freely the notations of \cite{HJ}, and ignore those from Sections~\ref{sec: 1}--\ref{sec: GT bases}. The quantum loop algebra $U_q(\Glie)$ admits Drinfeld generators $x_{i,r}^{\pm}$, $\phi_{i,\pm m}^{\pm}$ for $i \in I$, $r \in \BZ$ and $m \in \BZ_{\geq 0}$; see \cite[equation~(2.2)]{HJ}. ($\Glie$ here is $\widehat{\mathfrak{a}}$ in the introduction.) We shall need the following basic facts.
\begin{itemize}\itemsep=0pt
\item[(A)] Algebra $U_q(\Glie)$ is generated by $S := \{\phi_{i,\pm m}^{\pm}, x_{i,r}^+, x_{i,0}^- \, |\, i \in I, \, m \in \BZ_{\geq 0}, \, r \in \BZ \}$.
\item[(B)] For $i \in I$, $(x_{i,0}^+, x_{i,0}^-, \phi_0^+)$ generates $U_{q_i}(\mathfrak{sl}_2)$ with coproduct:
\begin{gather*} \Delta(x_{i,0}^+) = \phi_{i,0}^+ \otimes x_{i,0}^+ + x_{i,0}^+ \otimes 1,\qquad \Delta(x_{i,0}^-) = 1 \otimes x_{i,0}^- + x_{i,0}^- \otimes \phi_{i,0}^-,\\
\Delta(\phi_{i,0}^+) = \phi_{i,0}^+\otimes \phi_{i,0}^+. \end{gather*}
\end{itemize}

Let us be in the situation of \cite[Section~4.2]{HJ} where an inductive system of KR modules $(F_{k,l}\colon V_l \longrightarrow V_k)_{l<k}$ for a~f\/ixed $i \in I$ was constructed. For $j \in I$, $r \in \BZ$ and $m \in \BZ_{\geq 0}$, the $F_{k,l}$ commute with the $x_{j,r}^+$, and the $\phi_{j,\pm m}^{\pm}F_{k,l}$ for $k > l$, as Laurent polynomials in $q_i^k$ is described in \cite[Proposition~4.2]{HJ}. It is therefore enough to study the $x_{j,0}^- F_{k,l}$ for $k > l$.

If $j \neq i$, then $x_{j,0}^-$ annihilates the highest $\ell$-weight vector of $V_k$. This gives $x_{j,0}^- F_{k,l} = F_{k,l} x_{j,0}^-$. Assume $j = i$. The structural maps $F_{k,l}$ come from $U_q(\Glie)$-linear maps in \cite[Theorem~3.15]{HJ} which we write as $\SF_{k,l}\colon V_l \otimes Z_{lk} \longrightarrow V_k$, where $Z_{lk} := L(M_k^{\leq (-2l+1)d_i})$ before \cite[equation~(4.26)]{HJ} with highest $\ell$-weight vector $v_{lk}$ f\/ixed so that $v_l \otimes v_{lk} \mapsto v_k$.

\begin{Claim} Let $v \in V_l$. Then for $k > l+1$, we have $x_{i,0}^- F_{k,l}(v) \in F_{k,l+1}(V_{l+1})$ with
\begin{gather*}
x_{i,0}^- F_{k,l}(v) = F_{k,l+1}\left( q_i^{l-k} F_{l+1,l}(x_{i,0}^- v) + \frac{q_i^{k-l}-q_i^{l-k}}{q_i-q_i^{-1}} \SF_{l+1,l} (v \otimes x_{i,0}^- v_{l,l+1}) \right).
\end{gather*}
\end{Claim}

The proof is the same as that of Proposition~\ref{prop: first properties of inductive system}(3), based on~(B), and the fact that for $j \in I$ the weight space of $V_k$ of weight $\overline{\omega}_i^k \overline{\alpha}_j^{-1}$ is of dimension $\delta_{ji}$.

The asymptotic property of Section~\ref{sec: idea} is established for the generating set $S$ of $U_q(\Glie)$. One can then form a representation $\rho_c$ of $U_q(\Glie)$, with $c \in \BC^{\times}$ on the limit $V_{\infty}$ by specializing Laurent polynomials in $q_i^k$ to $c$. The representation $\rho_c$ contains a unique (up to multiple) vector annihilated by the $x_{j,r}^+$; it is $v_{\infty}$, the inductive limit of highest $\ell$-weight vectors $v_k \in V_k$. Moreover, we obtain from the proof of \cite[Proposition~4.2]{HJ} that
\begin{gather*} \phi_j^{\pm}(z) v_{\infty} = v_{\infty} \qquad \mathrm{for}\quad j \neq i,\qquad \phi_i^{\pm}(z) v_{\infty} = \frac{c-zac^{-1}}{1-za} v_{\infty}. \end{gather*}
The normalized $q$-character of $\rho_c$ is identical to that of $L_{i,a}^-$ in \cite[Theorem~6.1]{HJ}. The representation $\rho_c$ belongs to category $\BGG$ of $U_q(\Glie)$-modules \cite{He,Mukhin}.

In \cite[Lemma 4.4, Proposition~4.5]{HJ}, asymptotic property for the $x_{j,r}^-F_{k,l}$ was deduced from delicate results on $q$-characters of representations \cite{Hernandez}.

The arguments should work for Yangians: inductive systems of KR modules come from cyclicity result of particular tensor products of fundamental modules \cite{GT}; their asymptotic property is reduced eventually to representation theory of~$\mathfrak{sl}_2$. See \cite[Appendix~A]{FZ} for the case of Yangian of $\mathfrak{sl}_2$. Alternatively one may start from asymptotic modules $(\rho_c,V_{\infty})$ of a~quantum loop algebra and transform them into modules over the Yangian via the functor of Gautam--Toledano Laredo \cite[Section~6]{GTL}.

\subsection*{Acknowledgements}

The author thanks Vyjayanthi Chari, Giovanni Felder, David Hernandez, Masaki Kashiwara, Eugene Mukhin, Zengo Tsuboi, and Weiqiang Wang for interesting discussions. He is supported by the National Center of Competence in Research SwissMAP~-- The Mathematics of Physics of the Swiss National Science Foundation.

\pdfbookmark[1]{References}{ref}
\LastPageEnding

\end{document}